\documentstyle[amssymb,amsfonts]{amsart}

\def\R{{\hbox{\bf R}}}

\def\Q{{\hbox{\bf Q}}}

\def\P{{\hbox{\bf P}}}
\def\E{{\hbox{\bf E}}}

\font \roman = cmr10 at 10 true pt

\def\Vol{\hbox{\rm Vol}}

\def\be#1{ \begin{equation}\label{#1} }

\def\bas{\begin{align*}}
\def\eas{\end{align*}}
\def\bi{\begin{itemize}}
\def\ei{\end{itemize}}

\def\dim{{\hbox{\roman dim}}}
\def\rank{{\hbox{\roman rank}}}

\def\Z{{\hbox{\bf Z}}}

\def\eps{\varepsilon}
\newenvironment{proof}{\noindent {\bf Proof} }{\endprf\par}
\def \endprf{\hfill  {\vrule height6pt width6pt depth0pt}\medskip}
\def\emph#1{{\it #1}}
\def\textbf#1{{\bf #1}}

\def\Bv{{\mathbf v}}
\def\Bw{{\mathbf w}}

\def\BZ{{\mathbf Z}}

\def\BBP {{\mathbb P}}

\def\ep{{\epsilon}}

\theoremstyle{plain}
  \newtheorem{theorem}[subsection]{Theorem}

  \newtheorem{proposition}[subsection]{Proposition}
  
  \newtheorem{lemma}[subsection]{Lemma}
  \newtheorem{corollary}[subsection]{Corollary}

\theoremstyle{remark}
  \newtheorem{remark}[subsection]{Remark}
  
  \newtheorem{example}[subsection]{Example}

\theoremstyle{definition}
  \newtheorem{definition}[subsection]{Definition}

\include{psfig}

\begin{document}

\title[Inverse Littlewood-Offord and condition number]
{Inverse Littlewood-Offord theorems and the condition number of random discrete matrices}

\author{Terence Tao}
\address{Department of Mathematics, UCLA, Los Angeles CA 90095-1555}
\email{tao@@math.ucla.edu}
\thanks{T. Tao is a Clay Prize Fellow and is supported by a grant from the Packard Foundation.}

\author{Van H. Vu}
\address{Department of Mathematics, Rutgers, Piscataway, NJ 08854-8019}
\email{vanvu@@math.rutgers.edu}
\thanks{V. Vu is an A. Sloan Fellow and is supported by an NSF Career Grant.}

\begin{abstract}
Consider a random sum $\eta_1 v_1 + \ldots + \eta_n v_n$, where $\eta_1,\ldots,\eta_n$ are i.i.d. random signs and
$v_1,\ldots,v_n$ are integers.  The Littlewood-Offord problem asks to maximize
concentration probabilities such as $\P( \eta_1 v_1 + \ldots + \eta_n v_n = 0)$ subject
 to various hypotheses on the $v_1,\ldots,v_n$.  In this paper we develop an
 \emph{inverse} Littlewood-Offord theory (somewhat in the spirit of Freiman's inverse theory
  in additive combinatorics), which starts with the hypothesis that a concentration
  probability is large, and concludes that almost all of the $v_1,\ldots,v_n$ are
   efficiently contained in a generalized arithmetic progression.  As an application we give
   a new bound on the
magnitude  of the least singular value of a random Bernoulli matrix, which in turn provides
 upper tail estimates on the condition number.
\end{abstract}

\maketitle

\section{Introduction}

Let $\Bv$ be a multiset (allowing repetitions) of $n$ integers
$v_1,\dots, v_n$. Consider a class of discrete random walks
$Y_{\mu,\Bv}$ on the integers $\Z$, which start at the origin and
consist of $n$ steps, where at the $i^{th}$ step one moves
backwards or forwards with magnitude $v_i$ and probability
$\mu/2$, and stays at rest with probability $1-\mu$.  More
precisely:

\begin{definition}[Random walks]
For any $0 \leq \mu \leq 1$, let $\eta^{\mu} \in \{-1,0,1\}$
 denote a random variable which equals $0$ with probability $1-\mu$ and
$\pm 1$ with probability $\mu/2$ each.  In particular, $\eta^{1}$
is a random sign $\pm 1$, while $\eta^{0}$ is identically zero.

\noindent Given $\Bv$, we define $Y_{\mu,\Bv}$ to be the random variable
$$ Y_{\mu, \Bv} := \sum_{i=1}^n \eta_i^{\mu} v_i  $$
\noindent where the $\eta_i^{\mu}$ are i.i.d copies of
$\eta^{\mu}$.  Note that the exact enumeration $v_1,\ldots,v_n$ of the multiset is irrelevant.

\noindent The \emph{concentration probability} $\BBP_\mu(\Bv)$ of this random walk is defined to be the quantity
\begin{equation} \label{equa:defofPmu}  \BBP_\mu(\Bv) := \max_{a \in \Z} \P( Y_{\mu,\Bv} = a
).\end{equation}
Thus we have $0 < \BBP_\mu(\Bv) \leq 1$ for any $\mu, \Bv$.
\end{definition}

The concentration probability (and more generally, the
concentration function) is a central notion in probability theory
and has been studied extensively, especially by the Russian school
(see \cite{Sal, Rog, PJ} and the references therein).

The first goal of this paper is to establish a relation between
the magnitude of $\BBP_\mu(\Bv)$ and the arithmetic structure of
the multiset $\Bv=\{v_1, \dots, v_n\}$.   This gives an answer to
the general question of finding conditions under which one can squeeze
large probability inside a small interval. We will primarily be
interested in the case $\mu=1$, but for technical reasons it will
 be convenient to consider more general values of $\mu$.
Generally, however, we think of $\mu$ as fixed, while letting $n$
become very large.

A classical result of Littlewood-Offord  \cite{LO}, found in their
study of the number of real roots of random polynomials  asserts
that if all of the $v_i$'s are non-zero, then $\BBP_1(\Bv) =
O(n^{-1/2} \log n)$. The $\log$ term was later removed by
Erd\H{o}s \cite{Erd1}. Erd\H{o}s' bound is sharp, as shown by the
case $v_1=\dots=v_n \neq 0$. However, if one forbids this special
case and assumes that the $v_i$'s are all distinct, then the bound
can be improved significantly. Erd\H{o}s and Moser \cite{erdmos}
showed that under this stronger assumption, $\BBP_1(\Bv)
=O(n^{-3/2} \ln n)$. They conjectured that the logarithmic term is
not necessary and this was confirmed by S\'ark\"ozy and
Szemer\'edi \cite{SS}. Again, the bound is sharp (up to a constant
factor), as can be seen by taking $v_1,\ldots,v_n$ to be a proper
arithmetic progression such as $1,\ldots,n$. Later, Stanley
\cite{Stan}, using algebraic methods,  gave
  a very explicit bound for the probability in question.

The higher dimensional version of Littlewood-Offord's problem
(where the $v_i$ are non-zero vectors in $\R^d$, for some fixed
$d$) also drew lots of attention. Without the assumption that the
$v_i$'s are different, the best result was obtained by Frankl and
F\"uredi in \cite{FF}, following earlier results by Katona
\cite{Kat}, Kleitman \cite{Kle}, Griggs, Lagarias,  Odlyzko and
Shearer \cite{GLOS} and many others. However, the techniques used
in these papers did not seem to yield the generalization of
S\'ark\"ozy and Szemer\'edi's result (the $O(n^{-3/2})$ bound
under the assumption that the vectors are different).

 The generalization of S\'ark\"ozy and
Szemer\'edi's  result was  obtained by Hal\'asz \cite{Hal}, using analytical  methods (especially harmonic analysis). Hal\'asz' paper
was one of our starting points in this study.

\vskip2mm In the above two examples, we see that in order to make
$\BBP_{\mu}(\Bv)$ large, we have to impose a very strong additive
structure on $\Bv$ (in one case we set the $v_i$'s to be the same,
while in the other we set them to be elements of an arithmetic
progression). We are going to show that this is
the only way to make $\BBP_{\mu}(\Bv)$ large. More precisely, we
propose the following phenomenon:

\vskip2mm

 \centerline{\it If $\BBP_{\mu}(\Bv)$ is large, then $\Bv$ has a
strong additive structure. }

\vskip2mm In the next section, we are going to present several
theorems supporting this phenomenon. Let us mention here that
there is an analogous phenomenon in combinatorial number theory.
In particular, a famous theorem of Freiman asserts that if $A$ is
a finite set of integers and $A+A$ is small, then $A$ is contained
efficiently in a generalized arithmetic progression \cite[Chapter
5]{TVbook}. However, the proofs of Freiman theorem and those in
this paper are quite different.

 As an application, we are going to use these inverse theorems to study
random matrices. Let $M^{\mu}_n$ be an $n$ by $n$ random matrix,
whose entries are i.i.d copies of $\eta^{\mu}$. We are going to
show that with very high probability,
the condition number of $M^{\mu}_n$ is bounded from above by a polynomial in $n$
(see Theorem \ref{theo:condition2} below).  This result has high potential of applications in the theory of
probability in Banach spaces, as well as in numerical analysis and theoretical computer science.
A related result was recently established by Rudelson \cite{Rud}, with better upper bounds on the condition number but worse probabilities. We will discuss this application with more details  in Section 3.

To see the connection between this problem and inverse
Littlewood-Offord theory, observe that for any $\Bv=
(v_1,\ldots,v_n)$ (which we interpret as a column vector), the
entries of the product $M^\mu_n \Bv$ are independent copies of
$Y_{\mu,\Bv}$.  Thus we expect that $\Bv^T$ is unlikely to lie in
the kernel of $M^\mu_n$ unless the concentration probability
$\BBP_\mu(\Bv)$ is large.  These ideas are already enough to
control the singularity probability of $M^\mu_n$ (see e.g.
\cite{KKS, TV1, TV2}). To obtain the more quantitative condition
number estimates, we introduce a new discretization technique that
allows one to estimate the probability that a certain random
variable is small by the probability that a certain discretized
analogue of that variable is zero.

The rest of the paper is organized as follows. In Section \ref{sec2} we state our main
inverse theorems, and in Section \ref{sec3} we state our main results on condition numbers, as well as the key lemmas used to prove these results. In Section \ref{sec:appl}, we
give some brief applications of the inverse theorems.
In Section \ref{singsec} we prove the  result on condition numbers, assuming the inverse theorems and two other key ingredients:  a discretization
of generalized  progressions and an extension of the famous
result of Kahn, Koml\'os and Szemer\'edi \cite{KKS} on the probability that a
random Bernoulli matrix is singular.  The  inverse theorems are proven in
Section \ref{section:proofinverse}, after some preliminaries in Section
\ref{section:properties} in which we establish basic properties of $\BBP_{\mu}(\Bv)$. The result about  discretization
of progressions are proven  in Section \ref{discrete-sec}.
Finally in Section \ref{kks2-sec} we prove the extension of
 Kahn, Koml\'os and Szemer\'edi \cite{KKS}.

Let us conclude this section by setting out some basic notation. A
set $$P= \{c+ m_1a_1 + \dots + m_d a_d| M_i \le m_i \le M_i'\}$$
 is called a {\it generalized arithmetic
progression } (GAP) of rank $d$. It is convenient to think of $P$
as the image of an integer box $B:= \{(m_1, \dots, m_d)| M_i \le
m_i \le M_i' \}$ in $\Z^d$ under the linear map
$$\Phi: (m_1,\dots, m_d) \mapsto c+ m_1a_1 + \dots  + m_d a_d. $$ The
numbers $a_i$ are the {\it generators } of $P$. In this paper, all
GAPs have rational generators. A GAP is {\it proper}  if $\Phi$ is
one to one on $B$. The product $\prod_{i=1}^d (M_i'-M_i+1)$ is the
{\it volume} of $P$. If $M_i=-M'_i$ and $c=0$ (so $P=-P$) then we
say that $P$ is \emph{symmetric}.

For a set $A$ of reals and a positive integer $k$, we define the iterated sumset
$$kA := \{a_1+\dots + a_k |a_i \in A  \}.$$
One should take care to distinguish the sumset $kA$ from the dilate $k \cdot A$, defined for any real $k$ as
$$k \cdot A := \{ ka| a\in A \}. $$

We always assume that $n$ is sufficiently large. The asymptotic
notation $O()$, $o()$, $\Omega()$, $\Theta()$ is used under the assumption that $n \rightarrow \infty$.
Notation such as $O_d(f)$ means that the hidden constant in $O$
depends only on $d$.

We thank the referee for many detailed comments and corrections.

\section{Inverse Littlewood-Offord theorems}\label{sec2}

Let us start by presenting an example when $\BBP_{\mu}(\Bv)$ is
large. This example is the motivation of our inverse theorems.

\begin{example}\label{vex}
Let $P$ be a symmetric generalized arithmetic
progression of rank $d$ and volume $V$; we view $d$ as being fixed independently of $n$, though $V$ can grow with $n$.
Let $v_1,
\dots, v_n$ be (not necessarily different)  elements of $V$. Then
the random variable $Y_{\mu, \Bv} =\sum_{i=1}^n \eta_i v_i$ takes
values in the GAP $nP$ which has volume $n^d V$. From the pigeonhole principle it follows that
$$\BBP_{\mu} (\Bv) \ge n^{-d} V^{-1}.$$
In fact, the central limit theorem suggests that $\BBP_\mu(\Bv)$ should typically be of the order of
$n^{-d/2} V^{-1}$.
\end{example}

\noindent This example shows that if the elements of $\Bv$ belong
to a GAP with small rank and small volume then  $\BBP_{\mu} (\Bv)$ is
large. One might hope that the inverse  also holds, namely,

\vskip2mm

{\it If $\BBP_{\mu} (\Bv)$ is large, then (most of) the elements of $\Bv$
belong to a GAP with small rank and small volume. }

\vskip2mm In the rest of this section, we present three theorems,
which support this statement in a quantitative way.

\begin{definition}[Dissociativity]
Given a multiset $\Bw= \{w_1,\dots ,w_r\}$ of real numbers
and a positive number $k$,
we define the GAP $Q(\Bw,k)$ and the cube $S(\Bw)$ as follows:
\begin{align*}
Q(\Bw,k) &:= \{ m_1w_1 + \dots + m_r w_r| -k \le m_i \le k \} \\
S(\Bw) &:= \{\ep_1 w_1 +\dots + \ep_r w_r | \ep_i \in \{-1,1\} \}.
\end{align*}
We say that $\Bw$ is \emph{dissociated} if $S(\Bw)$ does not contain zero.
Furthermore,
 $\Bw$ is \emph{$k$-dissociated} if there do not exist integers $-k \leq m_1,\ldots,m_r \leq k$,
 not all zero, such that $m_1 w_1 + \ldots + m_r w_r = 0$.
\end{definition}

Our first result is the following simple proposition:

\begin{proposition}[Zeroth inverse theorem]\label{theo:inverse1}
 Let $\Bv = \{v_1,\ldots,v_n\}$ be
such that $\BBP_1(\Bv) > 2^{-d-1}$ for some integer $d \geq 0$. Then $\Bv$
contains a subset $\Bw$ of size $d$ such that the cube $S(\Bw)$
contains $v_1,\ldots,v_n$.
\end{proposition}

The next two theorems are more involved and also more useful.  In
these two theorems and their corollaries, we assume that $k$ and
$n$ are sufficiently large, whenever needed.

\begin{theorem}[First inverse theorem]\label{theo:inverse2}
Let $\mu$ be a positive constant at most 1 and  let $d$ be a
positive integer.  Then there is a constant $C=C(\mu,d) \geq 1$
such that the following holds. Let $k \geq 2$ be an integer and
let $\Bv = \{v_1, \ldots, v_n\}$ be a multiset such that
$$ \BBP_{\mu} ( \Bv )  \ge C(\mu,d) k^{-d} .$$
Then  there exists a $k$-dissociated multiset $\Bw = \{w_1,\ldots,w_r\}$ such that

\begin{enumerate}

\item $r \le d-1$ and $w_1,\ldots,w_r$ are elements of $\Bv$;

\item The union $\bigcup_{\tau \in \Z,  1 \leq \tau \leq k}
\frac{1}{\tau} \cdot Q(\Bw, k)$ contains all but $k^2$ of the
integers $v_1,\ldots,v_n$ (counting multiplicity). 

\end{enumerate}

\end{theorem}

This theorem should be compared against the heuristics in Example \ref{vex} (setting $k$ equal to
a small multiple of $\sqrt{n}$). In particular, notice
that the GAP $Q(\Bw, k)$ has very small volume, only $O(k^{d-1})$.

The above theorem does not yet show that most of the elements of
$\Bv$ belong to a single GAP. Instead, it shows that they belong
to the union of a few dilates of a GAP. One could remove the
unwanted $\frac{1}{\tau}$ factor by clearing denominators, but
this costs us an exponential factor such as $k!$, which is often
too large in applications. Fortunately, a more refined argument
allows us to eliminate these denominators while losing only
polynomial factors in $k$:

\begin{theorem}[Second inverse theorem]\label{theo:inverse3}
Let $\mu$ be a positive constant at most one, $\ep$ be an
arbitrary positive constant and $d$ be a positive integer. Then
there are constants $C=C(\mu,\ep,d) \geq 1$ and $k_0 = k_0 (\mu,
\ep, d) \ge 1$  such that the following holds. Let $k \geq k_0$ be
an integer and let $\Bv = \{v_1, \ldots, v_n\}$ be a multiset such
that

$$ \BBP_{\mu} ( \Bv )  \ge C k^{-d} .$$

  Then there exists a GAP
$Q$ with the following properties

\begin{enumerate}

 \item The rank of $Q$ is at most $d-1$;

 \item The  volume of $Q$ is
at most $ k^{2(d^2-1) +\ep}$;

\item $Q$ contains  all but at most $\ep k^2 \log k $ elements of
$\Bv$ (counting multiplicity);

\item There exists a positive integer $s$ at most $k^{d+\ep}$ such
that $su \in \Bv$ for each generator $u$ of $Q$. \end{enumerate}

\end{theorem}

\begin{remark} A
 small number of exceptional elements cannot be avoided.
 For instance, one can add $O(\log k)$ completely arbitrary elements to $\Bv$,
 and decrease $\P_\mu(\Bv)$ by a factor of $k^{-O(1)}$ at
 worst.
\end{remark}

For the applications in this paper, the following corollary of
Theorem \ref{theo:inverse3} is convenient.

\begin{corollary}\label{cor:inverse3} For any positive constants
 $A$ and $\alpha $ there is a positive constant
 $A'$ such that the following holds.
Let $\mu$ be a positive constant at most one and assume that $\Bv
= \{v_1, \ldots, v_n\}$ is a multiset of integers satisfying
$\BBP_{\mu} (\Bv) \geq n^{-A}$. Then there is a GAP $Q$ of rank at
most $A'$ and volume at most $n^{A'}$ which contains all but at
most $n^{\alpha}$ elements of $\Bv$ (counting multiplicity).
Furthermore, there exists a positive integer $s \le n^{A'}$ such
that $su \in \Bv$ for each generator $u$ of $\Q$.
\end{corollary}

\begin{remark}
The assumption $\BBP_{\mu} (\Bv) \geq n^{-A}$ in all statements
can be replaced by the following more technical, but somewhat
weaker assumption, that

$$\int_0^1 \prod_{i=1} |(1-\mu) + \mu \cos 2\pi v_i \xi | \,\, d \xi
\ge n^{-A}. $$

\noindent The right hand side is an upper bound for $\BBP_{\mu}
(\Bv)$, provided that $\mu$ is sufficiently small. Assuming that
$\BBP_{\mu} (\Bv) \ge n^{-A}$, what we will really use in the
proofs is the consequence

$$\int_0^1 \prod_{i=1} |(1-\mu) + \mu \cos 2\pi v_i \xi | d \xi
\ge n^{-A}. $$
 (See Section 5
for more details.) This weaker assumption is useful in
applications (see \cite{TVcond}).
\end{remark}

The vector versions of all three theorems (when the $v_i$'s are
vectors in $\R^r$, for any positive integer $r$) hold, thanks to
Freiman's isomorphism principle ( see, e.g., \cite[Chapter
5]{TVbook}). This principle allows us to project the problem from
$\R^r$ onto $\BZ$. The value of $r$ is irrelevant and does not
appear in any quantitative bound. In fact, one can even replace
$\R^r$ by any torsion free additive group.

Finally, let us mention that  in an earlier paper \cite{TV2} we
introduced another type of inverse Littlewood-Offord theorem. This
result showed that if $\BBP_\mu(\Bv)$ was comparable to
$\BBP_1(\Bv)$, then $\Bv$ could be efficiently contained inside a
GAP of bounded rank (see \cite[Theorem 5.2]{TV2} for details).

We shall prove these inverse theorems in Section \ref{section:proofinverse},
after some combinatorial and
Fourier-analytic preliminaries in Section \ref{section:properties}.
For now, let us take these results for granted and turn to an
application of these inverse theorems to random matrices.

\section{The condition number of random matrices}\label{sec3}

If $M$ is an $n \times n$ matrix, we use
$$\sigma_1 (M):= \sup_{x \in \R^n, \|x\|=1} \|Mx\|$$
to denote the largest singular value of $M$ (this parameter is also
often called the operator norm of $M$).  Here of course $\|x\|$ denotes the Euclidean magnitude of a vector $x \in \R^n$.
If $M$ is invertible, the \emph{condition number} $c(M)$ is defined as
$$c(M):= \sigma_1(M)  \sigma_1 (M^{-1}). $$
We adopt the convention that $c(M)$ is infinite if $M$ is not invertible.

The condition number  plays a
 crucial role in applied linear algebra and computer science. In
particular, the complexity of any algorithm which requires
solving a system of linear equations usually involves the
condition number of a matrix \cite{BT, Sma}. Another area of mathematics where
this parameter is important is the theory of probability in Banach spaces (see \cite{Lit, Rud}, for instance).

The condition number of a random matrix is a well-studied object (see \cite{Ede} and the references therein).
In the case when the entries of $M$ are i.i.d Gaussian random
variables (with mean zero and variance one), Edelman \cite{Ede}, answering a question of Smale \cite{Sma}
showed

\begin{theorem} \label{theo:condition1} Let $N_n$ be a $n \times n$ random matrix, whose
entries are i.i.d Gaussian random
variables (with mean zero and variance one). Then $\E(\ln c(N_n)) = \ln n + c+ o(1)$, where $c>0$ is an explicit constant.    \end{theorem}

In application, it is usually useful to have  a tail estimate. It was shown
by Edelman and Sutton  \cite{ET} that

\begin{theorem} \label{theo:condition11} Let $N_n$ be a $n$ by $n$ random matrix, whose
entries are i.i.d Gaussian random
variables (with mean zero and variance one). Then for any constant $A >0$,
$$\P( c(N_n) \ge n^{A+1}) = O_A( n^{-A}) . $$
\end{theorem}

On the other hand, for the other basic case when the entries are
i.i.d Bernoulli random variables (copies of $\eta^1$), the
situation is far from being  settled. Even to prove that the
condition number is finite with high probability is a non-trivial
task (see \cite{Kom1}). The techniques used to study  Gaussian
matrices rely heavily on the explicit joint distribution of the
eigenvalues. This distribution is not available for  discrete
models.

Using our inverse theorems, we can prove the following result, which is comparable to
Theorem \ref{theo:condition11}, and is another main result of this paper. Let
$M^{\mu}_n$ be the $n$ by $n$ random matrix whose entries are
i.i.d copies of $\eta^{\mu}$. In particular, the Bernoulli matrix
mentioned above is the case when $\mu=1$.

\begin{theorem} \label{theo:condition2} For any positive constant
$A$, there is a positive constant $B$ such that the following
holds. For any positive constant $\mu$ at most one and any
sufficiently large $n$

$$\P( c(M^{\mu}_n) \ge n^{B}) \le n^{-A}. $$
\end{theorem}

Given an invertible matrix $M$ of order
$n$, we set $\sigma_n (M)$ to be the smallest singular value of
$M$:
$$\sigma_n (M) := \min_{x \in \R^n, \|x\|=1} \|Mx\|. $$
Then we have
$$ c(M) = \sigma_1(M) / \sigma_n(M).$$

It is well known that there is a constant $C_\mu$ such that the
largest singular value  of $M^{\mu}_n$ is at most $C_\mu n^{1/2}$ with
exponential probability $1- \exp(- \Omega_\mu(n) )$ (see, for instance \cite{KV}). Thus, Theorem \ref{theo:condition2}
reduces to the following lower tail estimate for the smallest singular value of $\sigma_n(M)$:

\begin{theorem} \label{theo:singular1}
For any positive constant $A$, there is a positive constant $B$
such that the following holds. For any positive constant $\mu$ at
most one and any sufficiently large $n$

$$\P( \sigma_n(M^{\mu}_n) \le n^{-B}) \le n^{-A}.$$

\end{theorem}

Shortly prior to this paper, Rudelson \cite{Rud} proved the
following result.

\begin{theorem} \label{theo:singular2} Let $0 < \mu \leq 1$.
 There are positive constants $c_1(\mu), c_2(\mu)$ such that the following holds.
 For any $\ep \ge c_1(\mu) n^{-1/2}$
$$\P( \sigma_n(M^{\mu}_n) \le c_2(\mu) \ep n^{-3/2}) \le \ep. $$
\end{theorem}

In fact, Rudelson's result holds for a larger class of matrices.
The description of this class is, however, somewhat technical so
we refer the reader to \cite{Rud} for details.

It is useful to compare Theorems \ref{theo:singular1} and
\ref{theo:singular2}. Theorem \ref{theo:singular2} gives an
explicit dependence between the bound on $\sigma_n$ and the
probability, while the dependence between $A$ and $B$ in Theorem
\ref{theo:singular1} is implicit.  Actually our proof does provide
an explicit value for $B$, but it is rather large and we make no
attempt to optimize it. On the other hand, Theorem
\ref{theo:singular2} does not yield a probability better than
$n^{-1/2}$. In many applications (especially those involving the
union bound), it is important to have a probability bound of order
$n^{-A}$ with arbitrarily given $A$.






The proof of Theorem \ref{theo:singular1} relies on Corollary
\ref{cor:inverse3} and two other ingredients, which are of
independent interest. In the rest of this section, we discuss
these ingredients.  These ingredients will then be combined in Section \ref{singsec}
to prove Theorem \ref{theo:singular1}.

\subsection{Discretization of GAPs}\label{discrete-subsec}

Let $P$ be a GAP of integers of  rank $d$ and volume $V$. We show that given any
specified scale parameter $R_0$, one can ``discretize'' $P$ near
the scale $R_0$. More precisely, one can cover $P$ by the sum of a
coarse progression and a small progression, where the diameter of
the small progression is much smaller (by an arbitrarily specified
factor of $S$) than the spacing of the coarse progression, and
that both of these quantities are close to $R_0$ (up to a bounded
power of $SV$).

\begin{theorem}[Discretization]\label{discrete-thm}
 Let $P \subset \Z$ be a symmetric GAP of rank $d$ and volume $V$.  Let $R_0, S$ be positive integers.
Then there exists a scale $R \geq 1$ and two GAPs
$P_{\operatorname{small}}$, $P_{\operatorname{sparse}}$ of
\emph{rational} numbers with the following properties.
\begin{itemize}
\item(Scale)  $R = (SV)^{O_d(1)} R_0$. \item(Smallness)
$P_{\operatorname{small}}$ has rank at most $d$, volume at most
$V$, and takes values in $[-R/S, R/S]$. \item(Sparseness)
$P_{\operatorname{sparse}}$ has rank at most $d$, volume at most
$V$, and any two distinct elements of $SP_{\operatorname{sparse}}$
are separated by at least $RS$. \item(Covering)  $P \subseteq
P_{\operatorname{small}} + P_{\operatorname{sparse}}$.
\end{itemize}
\end{theorem}

This theorem is elementary but is somewhat involved and the
detailed proof will appear in Section \ref{discrete-sec}. Let us,
at this point, give an informal explanation,  appealing to the
analogy between the combinatorics of progressions and linear
algebra. Recall that a  GAP of rank $d$ is the image $\Phi( B )$
of a $d$-dimensional box under a linear map $\Phi$. This can be
viewed as a discretized, localized analogue of the object
$\Phi(V)$, where $\Phi$ is a linear map from a $d$-dimensional
vector space $V$ to some other vector space.  The analogue of a
``small'' progression would be an object $\Phi(V)$ in which $\Phi$
vanished. The analogue of a ``sparse'' progression would be an
object $\Phi(V)$ in which the map $\Phi$ was injective.  Theorem
\ref{discrete-thm} is then a discretized, localized analogue of
the obvious linear algebra fact that given any object of the form
$\Phi(V)$, one can split $V = V_{\operatorname{small}} +
V_{\operatorname{sparse}}$ for which
$\Phi(V_{\operatorname{small}})$ is small and
$\Phi(V_{\operatorname{sparse}})$ is sparse.  Indeed one simply
sets $V_{\operatorname{small}}$ to be the kernel of $\Phi$, and
$V_{\operatorname{sparse}}$ to be any complementary subspace to
$V_{\operatorname{small}}$ in $V$.  The proof of Theorem
\ref{discrete-thm} that we give follows these broad ideas, with
$P_{\operatorname{small}}$ being essentially a ``kernel'' of the
progression $P$, and $P_{\operatorname{sparse}}$ being a kind of
``complementary progression'' to this kernel.

To oversimplify enormously, we shall exploit this discretization result
(as well as the inverse Littlewood-Offord theorems) to control the event that the singular value is small, by
the event that the singular value (of a slightly modified random matrix)
is \emph{zero}.  The control of this latter quantity is the other ingredient
of the proof, to which we now turn.

\vskip2mm

\subsection{Singularity of random matrices}

A famous result of
Kahn, Koml\'os and Szemer\'edi \cite{KKS} asserts that the
probability that $M^{1}_n$ is singular (or equivalently, that $\sigma_n(M^1_n) = 0$)
is exponentially small:

\begin{theorem} \label{theo:KKS1} There is a positive constant
$\eps$ such that

$$\P( \sigma_n(M^1_n) = 0 ) \le (1-\eps )^n. $$
\end{theorem}

In \cite{KKS} it was shown that one can take $\eps=.001$.
Improvements on $\eps$ are obtained recently in \cite{TV1, TV2}.
The value of $\ep$ does not play a critical role in this paper.

 To prove Theorem \ref{theo:condition2}, we need the following
generalization of Theorem \ref{theo:KKS1}. Notice that the row
vectors of $M^1_n$ are i.i.d copies of $X^1$, where
$X^1=(\eta^1_1, \dots, \eta^1_n)$ and $\eta_i^1$ are i.i.d copies
of $\eta^1$. By changing $1$ to $\mu$, we can define $X^{\mu}$ in
the obvious manner.  Now let $Y$ be a set of $l$ vectors $y_1,
\dots, y_l$ in $\R^n$ and $M^{\mu, Y}_n$ be the random matrix
whose rows are $X_1^{\mu}, \dots, X_{n-l}^{\mu}, y_1, \dots, y_l$,
where $X_i^{\mu}$ are i.i.d copies of $X^{\mu}$.

\begin{theorem} \label{theo:KKS2} Let $0 < \mu \leq 1$, and let $l$ be a
non-negative  integer. Then there is a positive constant
$\eps=\eps(\mu, l)$ such that the following holds. For any set $Y$
of $l$ independent vectors from $\R^n$,
$$\P( \sigma_n(M^{\mu,Y}_n) = 0 ) \le (1-\eps)^n. $$

\end{theorem}

\begin{corollary} \label{cor:KKS2} Let $0 < \mu \leq 1$.
Then there is a positive constant $\eps=\eps(\mu)$ such that the
following holds. For any vector $y\in \R^n$, the probability that
there are $w_1, \dots, w_{n-1}$, not all zeros, such that
$$y = X_1^{\mu} w_1 + \dots X_{n-1}^{\mu} w_{n-1} $$
\noindent is at most $(1-\eps)^n$.
\end{corollary}

We will prove Theorem \ref{theo:KKS2} in Section \ref{kks2-sec} by using the machinery from \cite{TV1}.

\section{Some quick applications of the inverse theorems}\label{sec:appl}

The inverse theorems provide effective bounds for counting the number of ``exceptional'' collections $\Bv$ of numbers with high concentration probability; see for instance \cite{TV2} for a demonstration of how such bounds can be used in applications.
In this section,  we present two such bounds that can be obtained from the inverse theorems developed here. In the first
example, let $\ep$ be a positive constant and $M$ be a large
integer and consider the following question:

{\it How many sets $\Bv$ of $n$ integers with absolute values at
most $M$ are there such that $\BBP_1(\Bv) \ge \ep $ ? }

By Erd\H{o}s' result, all but at most $O(\ep^{-2})$ of the
elements of $\Bv$ are non-zero. Thus we have the upper bound ${n
\choose {\ep^{-2}} } (2M+1)^{O(\ep^{-2})}$ for the number in
question. Using Proposition \ref{theo:inverse1}, we can obtain a
better bound as follows. There are only $M^{O(\ln \ep^{-1})}$ ways
to choose the generators of the cube. After the cube is fixed, we
need to choose $O(\ep^{-2})$ non-zero elements inside it. As the
cube has volume $O(\ep^{-1})$, the number of ways to do this is
$(\frac{1}{\ep})^{O(\ep^{-2})}$. Thus, we end up with a bound

$$ M^{O(\ln \ep^{-1})} (\frac{1}{\ep})^{O(\ep^{-2})} $$

\noindent which is better than the previous one if $M$ is
considerably larger than $\ep^{-1}$.

For the second application, we return to the question of bounding the singularity
probability $\P( \sigma_n(M^1_n) = 0 )$ studied in Theorem \ref{theo:KKS1}.  This probability
is conjectured to equal $(1/2+o(1))^n$, but this remains open (see \cite{TV2} for the
latest results and some further discussion). The event that $M^1_n$ is singular is the
same as the event that there exists some non-zero vector $v \in
\R^n$ such that $M^1_n v=0$. For simplicity, we use the notation
$M_n$ instead of $M^1_n$ in the rest of this section. It turns out
that one can obtain the optimal bound $(1/2+o(1))^n$ if one
restricts $v$ to some special set of vectors.

Let $\Omega_1$ be the set of vectors in $\R^n$ with at least $3n
/\log_2 n$ coordinates. Koml\'os proved the
following:

\begin{theorem}\label{theo:komlos-transpose}
The probability  that $M_n v=0$ for  some non-zero $v\in \Omega_1$
is  $(1/2+o(1))^n$.
\end{theorem}

A proof of this theorem can be found in Bollob\'as' book \cite{Bol}.

We are going to consider another restricted class. Let $C$ be an
arbitrary positive  constant and let  $\Omega_2$ be the set of
integer vectors in $\R^n$ where the coordinates have absolute
values at most $n^C$. Using Theorem \ref{theo:inverse2}, we can
prove

\begin{theorem} \label{theo:Omega2} The probability that $M_n v=0$
for some non-zero $v \in \Omega_2$ is $(1/2+o(1))^n$.
\end{theorem}

\begin{proof} The lower bound is trivial so we focus on the upper
bound. For each non-zero vector $v$, let $p(v)$ be the probability
that $X \cdot v=0$, where $X$ is a random Bernoulli vector. From
independence we have $\P (M_n v =0) =p(v)^n$. Since a hyperplane
can contain at most $2^{n-1}$ vectors from $\{-1,+1\}^n$, $p(v)$
is at most $1/2$. For $j=1,2, \dots$, let $S_j$ be the number of
non-zero vectors $v$ in $\Omega_2$ such that $2^{-j-1} < p(v) \le
2^{-j}$. Then the probability that $M_n v=0$ for some non-zero $v
\in \Omega_2$ is at most

$$\sum_{j=1}^{n} (2^{-j})^n S_j. $$

\noindent Let us now restrict the range of $j$. Notice that if
$p(v) \ge n^{-1/3}$, then by Erd\H{o}s's result (mentioned in the
Introduction)  most of the coordinates of $v$ are zero.  In this
case, by Theorem \ref{theo:komlos-transpose} the contribution from
these $v$ is at most $(1/2+o(1))^n$. Next, since the number of
vectors in $\Omega_2$ is at most $(2n^C +1)^n \le n^{(C+1)n}$, we
can ignore those $j$ where $2^{-j} \le n^{-C-2}$. Now it suffices
to show

$$\sum_{n^{-C-2} \le 2^{-j} \le n^{-1/3}} (2^{-j})^n S_j =o( (1/2)^n).
$$

 For any
relevant $j$, we can find an integer $d=O(1)$ and a positive
number $\ep =\Omega (1)$ such that
$$ n^{-(d - 1/3) \ep}  \le 2^{-j} < n^{-(d -2/3)\ep}. $$

 Set
$k := n^{\ep}$. Thus $2^{-j} \gg k^{-d}$ and we can use Theorem
\ref{theo:inverse2} to estimate $S_j$. Indeed, by invoking this
theorem, we see that there are at most ${n \choose k^2}
(2n^{C}+1)^{k^2} = n^{O(k^2}) =n^{o(n)}$ ways to choose the
positions and values of exceptional coordinates of $v$.
Furthermore, There is only $(2n^{C}+1)^{d-1} = n^{O(1)}$ ways to
fix the generalized progression $P:= Q(\Bw, k)$.

Notice that the elements of $P$ are polynomially bounded in $n$.
Such integers have only $n^{o(1})$ divisors. So if $P$ is fixed
then  any (non-exceptional) coordinate of $v$ has at most $|P|
n^{o(1)}$ possible values. This means that once $P$ is fixed, the
number of ways to set the non-exceptional   coordinates of $v$ is
at most $(n^{o(1)} |P|)^n = (2k+1)^{(d-1 +o(1))n}$. Putting these
together,

$$ S_j \le  n^{O(k^2)} k^{(d-1 +o(1)) n}. $$

\noindent As $k=n^{\eps}$ and $2^{-j} \le n^{-(d- 2/3)\ep}$, it
follows that

$$ 2^{-jn} S_j \le   n^{o(n)} n^{- \ep n/3} = o(\frac{1}{\log n})2^{-n}.
$$

\noindent Since there are only $O(\log n)$ relevant $j$, we can
conclude the proof by summing the bound over $j$.
\end{proof}

\section{Properties of $\BBP_{\mu}(\Bv)$} \label{section:properties}

In order to prove the inverse Littlewood-Offord theorems in Section \ref{sec2}, we shall first
need to develop some useful tools for estimating the quantity $\BBP_\mu(\Bv)$.  That shall be the purpose of
this section.  We remark that the tools here are only used for the proof of the inverse Littlewood-Offord theorems in
Section \ref{section:proofinverse} and are not required elsewhere in the paper.

It is convenient to think of $\Bv$ as a word, obtained by
concatenating the numbers $v_i$:
$$\Bv = v_1 v_2 \dots v_n. $$
This will allow us to perform several operations such as
concatenating, truncating and repeating. For instance, if $\Bv
=v_1\dots v_n$ and $\Bw= w_1\dots w_m$, then
$$\BBP_{\mu} (\Bv \Bw) = \max _{a \in Z} \Big(\sum_{i=1}^n \eta_i^{\mu} v_i
+ \sum_{j=1}^m \eta_{n+j} ^{\mu} w_j=a \Big)$$ \noindent where
$\eta_{k}^{\mu}, 1\le k \le n+m$ are  i.i.d copies of
$\eta^{\mu}$. Furthermore, we use $\Bv^k$ to denote the
concatenation of $k$ copies of $\Bv$.

It turns out that there is a nice calculus concerning the
expressions $\BBP_\mu(\Bv)$, especially when $\mu$ is small. The core
properties are summarized in the next lemma.

\begin{lemma} \label{lemma:properties} The following properties
hold.

\begin{itemize}

\item $\BBP_\mu(\Bv)$ is  invariant under permutations of $\Bv$.

\item For any words  $\Bv, \Bw$
\begin{equation}\label{pmub}
\BBP_\mu(\Bv) \BBP_\mu(\Bw) \leq \BBP_\mu(\Bv \Bw) \leq \BBP_\mu(\Bv).
\end{equation}

\item For any $0 < \mu \le 1$, any $0 < \mu' \leq \mu/4$, and any word $\Bv$,
\begin{equation}\label{14} \BBP_{\mu}(\Bv) \leq \BBP_{\mu'}(\Bv).
\end{equation}

\item  For any number $0 < \mu \leq 1/2$ and any word $\Bv$,
\begin{equation}\label{pmub-2} \BBP_\mu(\Bv) \leq \BBP_{\mu/k}(\Bv^k).
\end{equation}

\item For any number $0 < \mu \leq 1/2$ and any words $\Bv, \Bw_1,
\dots, \Bw_m$ we have
\begin{equation}\label{pmub-3}
\BBP_\mu(\Bv \Bw_1 \ldots \Bw_m) \leq \left( \prod_{j=1}^m \BBP_\mu(\Bv
\Bw_j^{m}) \right)^{1/m}.
\end{equation}

\item  For any number $0 < \mu \leq 1/2$ and any words $\Bv,
\Bw_1, \dots, \Bw_m$, there is an index $1\le j \le m$ such that
\begin{equation}\label{pmub-3a}
\BBP_\mu(\Bv \Bw_1 \ldots \Bw_m) \leq \BBP_\mu(\Bv \Bw_j^{m}).
\end{equation}

\end{itemize}

\end{lemma}

\begin{proof} The first two properties are trivial. To verify the
rest,  let us notice from Fourier analysis that

\begin{equation}\label{emu}
 \P( \eta^{(\mu)}_1 v_1 + \ldots + \eta^{(\mu)}_n v_n = a ) =
 \int_0^1 e^{-2\pi i a \xi} \prod_{j=1}^n (1-\mu + \mu \cos(2\pi v_j \xi))\ d \xi.
 \end{equation}

When $0 < \mu \leq 1/2$, the expression $1-\mu + \mu \cos(2\pi v_j
\xi))$  is positive, and we thus have

\begin{equation}\label{muf}
\BBP_\mu(\Bv) =\P( Y_{\mu,\Bv}=0) = \int_0^1 \prod_{j=1}^n (1-\mu +
\mu \cos(2\pi v_j \xi)) \ d \xi.
\end{equation}

To prove  \eqref{14}, notice that for any $0 < \mu \le 1$, $0 <
\mu' \leq \mu/4$ and any $\theta$ we have the elementary
inequality
$$|(1-\mu) + \mu \cos \theta| \le (1-\mu') +
\mu' \cos 2 \theta. $$

\noindent Using this, we have

\begin{align*}
\BBP_\mu(\Bv) &\le  \int_0^1 \prod_{j=1}^n |(1-\mu + \mu \cos(2\pi
v_j \xi))|  \ d \xi \\ &\le \int_0^1 \prod_{j=1}^n
(1-\mu' + \mu' \cos(4\pi v_j \xi))  \ d \xi \\
&= \int_0^1 \prod_{j=1}^n (1-\mu' + \mu'
\cos(4\pi v_j \xi)) d \xi\\
&= \BBP_{\mu'} (\Bv)
\end{align*}
\noindent where the next to last equality follows by changing
$\xi$ to $2\xi$ and considering the periodicity of cosine.

Similarly, observe that for $0 < \mu \leq 1/2$ and
$k \geq 1$ we have
$$ (1 - \mu + \mu \cos(2\pi v_j \xi)) \leq (1 - \frac{\mu}{k} + \frac{\mu}{k} \cos(2\pi v_j \xi))^k.$$

\noindent Indeed from the concavity of $\log(1-t)$ when $0 < t <
1$, we have $\log(1-t) \leq k \log(1-\frac{t}{k})$, and the claim
follows by exponentiating this with $t := \mu(1 - \cos(2\pi v_j
\xi))$). This proves \eqref{pmub-2}.

Finally, \eqref{pmub-3} is a consequence of \eqref{muf} and H\"older's
inequality, while \eqref{pmub-3a} follows directly from \eqref{pmub-3}.
\end{proof}

Now we consider the distribution of the equal-steps random walk
$\eta^{\mu}_1+ \dots +\eta_{m}^{\mu} = Y_{\mu,1^m}$. Intuitively,
this random walk is concentrated in an interval of length $O((1
+\mu m)^{1/2})$ and has a roughly uniform distribution in the
integers in this interval (though when $\mu$ is close to $1$,
parity considerations may cause $Y_{\mu,1^m}$ to favor the even
integers over the odd ones, or vice versa); compare with the
discussion in Example \ref{vex}. The following lemma is a
quantitative version of this intuition.

\begin{lemma}\label{lemma:torsion}  For any $0 < \mu \leq 1$ and $m \geq 1$ we have
\begin{equation}\label{pmu}
\BBP_{\mu}(1^m) = \sup_a \P( \eta^\mu_1 + \ldots + \eta^\mu_m = a
) = O( (\mu m)^{-1/2} ).
\end{equation}
In fact, we have the more general estimate
\begin{equation}\label{pma}
 \P( \eta^\mu_1 + \ldots + \eta^\mu_m = a ) = O( (\tau^{-1} + (\mu m)^{-1/2})
\P( \eta^\mu_1 + \ldots + \eta^\mu_m \in  [a -\tau, a+\tau] )
\end{equation}
for any $a \in \Z$ and $\tau \geq 1$.

Finally, if $\tau \geq 1$ and $S$ is any $\tau$-separated set of
integers (i.e. any two distinct elements of $S$ are at least
$\tau$ apart) then
\begin{equation}\label{pmus} \P( \eta^\mu_1 + \ldots + \eta^\mu_m \in S ) \leq O( \tau^{-1} + (\mu m)^{-1/2} ).
\end{equation}
\end{lemma}

\begin{proof}  We first prove \eqref{pmu}.
From \eqref{14} we may assume $\mu \leq 1/4$, and then by
\eqref{muf} we have
$$ \BBP_{\mu}(1^m) = \int_0^1 |1-\mu + \mu \cos(2\pi \xi)|^m\ d\xi.$$
Next we  use the elementary estimate
$$1-\mu + \mu \cos(2\pi \xi) \leq \exp( - \mu \|\xi\|^2/100 ),$$
where $\|\xi\|$ denotes the distance to the nearest integer. This
implies that $ \BBP_{\mu}(1^m)$ is bounded from above by $\int_0^1
\exp( - \mu m \|\xi\|^2/100 ) d \xi $, which is of order $O( (\mu
m)^{-1/2})$ (to see this notice that for $\xi \ge 1000(\mu
m)^{-1/2}$ the function $\exp( - \mu m \|\xi\|^2/100)$ is quite
small and its integral is negligible).

Now we prove \eqref{pma}.  We may assume that $\tau \leq (\mu
m)^{1/2}$, since the claim for larger $\tau$ follows
automatically. By symmetry we can take $a \geq 2$.

For each integer $a$, let $c_a$ denote the probability
$$ c_a := \P (\eta^{(\mu)}_1 + \ldots + \eta^{(\mu)}_m = a).$$
Direct computation (letting $i$ denote the number of
$\eta^{(\mu)}$ variables which equal zero) yields the explicit
formula
$$ c_a = \sum_{j=0}^m \binom{m}{j} (1-\mu)^j (\mu/2)^{m-j}
\binom{m-j}{(a+m-j)/2},$$ with the convention that the binomial
coefficient $\binom{a}{b}$ is zero when $b$ is not an integer
between $0$ and $a$. This in particular yields the monotonicity
property $c_a \geq c_{a+2}$ whenever $a \geq 0$. This is already
enough to yield the claim when $a > \tau$, so it remains to verify
the claim when $a \le \tau$.  Now the random variable $\eta^\mu_1
+ \ldots + \eta^\mu_m$ is symmetric around the origin and has
variance $\mu m$, so from Chebyshev's inequality we know that

$$ \sum_{0 \leq a \leq 2(\mu m)^{1/2}} c_a
=\Theta( 1).$$

From \eqref{pmu} we also have $c_a = O( (\mu m)^{-1/2} )$ for all
$a$.  From this and the monotonicity property $c_a \geq c_{a+2}$
and the pigeonhole principle we see that $c_a =\Theta( (\mu
m)^{-1/2})$ either for all even $0 \leq a \leq ( \mu m )^{1/2}$,
or for all odd $0 \leq a \leq (\mu m)^{1/2}$.  In either case, the
claim \eqref{pma} is easily verified.  The bound in \eqref{pmus}
then follows by summing \eqref{pma} over all $a \in S$ and noting
that $\sum_a c_a = 1$.
\end{proof}

One can also use the formula for $c_a$ to prove \eqref{pmu} as
well. The simple details are left as an exercise.

\section{Proofs of the inverse theorems }
\label{section:proofinverse}

We now have enough machinery to prove the inverse Littlewood-Offord theorems.  We first give a quick proof
of Proposition \ref{theo:inverse1}:

\begin{proof}[of Proposition \ref{theo:inverse1}]
Suppose that the conclusion failed.  Then an easy greedy algorithm argument shows that
 $\Bv$ must contain a dissociated subword
$\Bw = (w_1,\ldots,w_{d+1})$ of length $d+1$.  By \eqref{pmub}, we have
$$ 2^{-d-1} <  \BBP_1 (\Bv) \le \BBP_1(\Bw ).$$
On the other hand, since $\Bw$ is dissociated, all the sums of the
form $\eta_1 w_1 + \dots \eta_{d+1} w_{d+1} $ are distinct and so
$\BBP_1 (\Bw) \leq 2^{-d-1}$,  yielding the desired
contradiction.
\end{proof}

To prove Theorem \ref{theo:inverse2}, we modify the above argument by replacing the notion of dissociativity
by $k$-dissociativity.  Unfortunately this makes the proof somewhat longer:

\begin{proof}[of Theorem \ref{theo:inverse2}]
We construct an $k$-dissociated tuple $(w_1, \ldots, w_r)$ for some
$0 \leq r \leq d-1$ by the following algorithm:

\begin{itemize}

\item Step 0.  Initialize $r = 0$.  In particular, $(w_1, \ldots,
w_r)$ is trivially $k$-dissociated. From \eqref{pmub-2} we have
\begin{equation}\label{nbig-recurse} \BBP_{\mu/4d} (\Bv^d) \geq \BBP_{\mu/4} (\Bv ) \geq \BBP_\mu(\Bv).
\end{equation}

\item Step 1.   Count how many $1 \leq j \leq n$ there are such
that $(w_1, \ldots, w_r, v_j)$ is $k$-dissociated. If this number
is less than $k^2$, halt the algorithm.  Otherwise, move on to
Step 2.

\item Step 2.   Applying the last property of Lemma
\ref{lemma:properties}, we can locate a $v_j$ such that $(w_1,
\ldots, w_r, v_j)$ is $k$-dissociated, and
\begin{equation}\label{nbig-recurse2}
\BBP_{\mu/4d} (\Bv^{d-r} w_1^{k^2} \ldots w_r^{k^2}) \leq
\BBP_{\mu/4d} (\Bv^{d-r-1} w_1^{k^2} \ldots w_r^{k^2} v_j^{k^2}).
\end{equation}

We then set $w_{r+1} := v_j$ and increase $r$ to $r+1$. Return to
Step 1. Note that $(w_1, \ldots, w_r)$ remains $k$-dissociated,
and \eqref{nbig-recurse} remains true.

\end{itemize}

 Suppose  that we terminate at some step $r \le d-1$. Then we have an
$r$-tuple $(w_1,\ldots,w_r)$ which is $k$-dissociated, but such
that $(w_1,\ldots,w_r,v_j)$ is $k$-dissociated for at most $k^2$
values of $v_j$.  Unwinding the definitions, this shows that for
all but at most $k^2$ values of $v_j$, there exists $\tau \in
[1,k]$ such that $\tau v_j \in Q(\Bw, k)$, proving the claim.

It remains to show that  we must indeed terminate at some step $r
\le d-1$. Assume (for a contradiction) that we have reached step
$d$. Then we have an $k$-dissociated tuple $(w_1,\ldots,w_d)$, and by \eqref{nbig-recurse},
\eqref{nbig-recurse2} we have
$$ \BBP_{\mu} (\Bv ) \leq \BBP_{\mu/4d} ({w_1^{k^2} \ldots w_d^{k^2}}) = \P( Y_{\mu/4d, w_1^{k^2} \ldots w_d^{k^2} } = 0 ).$$
Let $\Gamma \subset \Z^d$ be the lattice
$$ \Gamma := \{ (m_1,\ldots,m_d) \in \Z^d: m_1 w_1 + \ldots + m_d w_d = 0 \},$$
then by using independence we can write
\begin{equation}\label{xid}
 \BBP_{\mu}(\Bv ) \leq \P( Y_{\mu/4d, w_1^{k^2} \ldots w_d^{k^2} } = 0 ) =
 \sum_{(m_1,\ldots,m_d) \in \Gamma} \prod_{j=1}^d \P( Y_{\mu/4d,1^{k^2}} = m_j
 ).\end{equation}

Now we use a volume packing argument. From Lemma \ref{lemma:torsion} we have
$$ \P( Y_{\mu/4d, 1^{k^2}} = m ) = O_{\mu,d}( \frac{1}{k} \sum_{m' \in m + (-k/2,k/2)}
 \P( Y_{\mu/4d, 1^{k^2}} = m' ) )$$
and hence from \eqref{xid} we have
\begin{align*}
\BBP_{\mu} (\Bv)  &\leq O_{\mu,d}( k^{-d}
\sum_{(m_1,\ldots,m_d) \in \Gamma} \\
&\quad \sum_{(m'_1,\ldots,m'_d) \in (m_1,\ldots,m_d) +
(-k/2,k/2)^d} \prod_{j=1}^d \P( Y_{\mu/4d, 1^{k^2}} = m_j' ) ).
\end{align*}
Since $(w_1,\ldots,w_d)$ is $k$-dissociated, all the $(m_1',\dots,
m_d')$ tuples in $\Gamma + (-k/2,k/2)^d$ are different. Thus, we
conclude
$$  \BBP_{\mu}(\Bv) \leq O_{\mu, d}\Big( k^{-d}
\sum_{(m_1,\ldots,m_d) \in \Z^d} \prod_{j=1}^d \P( Y_{\mu/4d,
1^{k^2}} = m_j ) \Big).$$ But from the union bound we have
$$ \sum_{(m_1,\ldots,m_d) \in \Z^d} \prod_{j=1}^d \P( Y_{\mu/4d, 1^{k^2}} = m_j ) =
1,$$
so
$$   \BBP_{\mu}(\Bv) \leq O_{\mu, d}( k^{-d}). $$
\noindent To complete the proof, set the constant $C = C(\mu,d)$ in the
theorem to be larger than the hidden constant in $O_{\mu,d}(k^{-d})$.
\end{proof}

\begin{remark} \label{remark:chernoff} One can also use the Chernoff bound and obtain a shorter proof (avoiding the volume packing argument) but with an extra logarithmic loss in the estimates.
\end{remark}

Finally we perform some additional arguments to eliminate the $\frac{1}{\tau}$ dilations in Theorem \ref{theo:inverse2} and obtain
our final inverse Littlewood-Offord theorem.  The key will be the following lemma.

Given a set $S$ and a number $v$. The torsion of $v$ with respect
to $S$ is the smallest positive integer $\tau$ such that $\tau v
\in S$. If such $\tau$ does not exists, we say that $v$ has
infinite torsion with respect to $S$.

 The key new ingredient will
be the following lemma, which asserts that adding a high torsion
element to a random walk  reduces the concentration probability
significantly.

\begin{lemma}[Torsion implies dispersion]\label{tor}
Let $0 < \mu \leq 1$ and consider a GAP $Q:= \{\sum_{i=1}^d x_i
W_i| -L_i \le x_i \le L_i \}$.
 Assume that $W_{d+1}$ has finite torsion $\tau$ with respect to $2Q$.
  Then there is a constant $C_{\mu}$ depending only on $\mu$ such
that
$$ \BBP_\mu( W_1^{L_1} \ldots W_d^{L_d} W_{d+1}^{\tau^2} )
\le C_{\mu} \tau^{-1} \BBP_\mu( W_1^{L_1} \ldots W_d^{L_d} ).$$
\end{lemma}

\begin{proof} Let $a$ be an integer such that

$$\BBP_\mu( W_1^{L_1} \ldots W_d^{L_d} W_{d+1}^{\tau^2} ) =
\P( \sum_{i=1}^{d} W_i \sum_{j=1}^{L_i} \eta^\mu_{j,i} +
 W_{d+1} \sum_{j=1}^{\tau^2} \eta^\mu_{j,d+1}   = a ), $$
 where the $\eta^\mu_{j,i}$ are i.i.d. copies of $\eta^\mu$.
 It suffices to show that

$$ \P( \sum_{i=1}^{d} W_i \sum_{j=1}^{L_i} \eta^\mu_{j,i} +
 W_{d+1} \sum_{j=1}^{\tau^2} \eta^\mu_{j,d+1}   = a )
 = O_{\mu}(\tau^{-1}) \BBP_\mu( W_1^{L_1} \ldots W_{d}^{L_d}  ).$$

 Let $S$ be the set of all $m \in [-\tau^2,\tau^2]$ such that
$Q + m W_{d+1}$ contains $a$.  Observe that in order for
$\sum_{i=1}^{d} W_i \sum_{j=1}^{L_i} \eta^\mu_{j,i} +
 W_{d+1} \sum_{j=1}^{\tau^2} \eta^\mu_{j,d+1} $ to
equal $a$, the quantity $\sum_{j=1}^k \eta^{\mu}_{j,d+1}$ must lie
in $S$. By  the definition of $\BBP_\mu( W_1^{L_1} \ldots
W_{d}^{L_d} )$ and Bayes identity, we conclude

$$  \P( \sum_{i=1}^{d} W_i \sum_{j=1}^{L_i} \eta^\mu_{j,i} +
 W_{d+1} \sum_{j=1}^{\tau^2} \eta^\mu_{j,d+1}   = a ) \leq  \BBP_\mu( W_1^{L_1} \ldots W_{d}^{L_d} )
\P( \sum_{j=1}^{\tau^2}  \eta^\mu_{j,d+1} \in S ).$$

Consider two elements $x,y \in S$. By the definition of $S$,
$(x-y)v \in Q-Q= 2Q$. From definition of $\tau$, $|x-y|$ is either
zero or at least $\tau$.  This implies that  $S$ is
$\tau$-separated and the claim now follows from Lemma
\ref{lemma:torsion}.
\end{proof}

We will also need the following technical lemma.

\begin{lemma} \label{lemma:sumofQ}
Consider a GAP $Q(\Bw, L)$. Assume that $v$ is an element with
(finite) torsion $\tau$ with respect to $Q(\Bw, L)$. Then

$$Q(\Bw, L) + Q(v, L') \subset \frac{1}{\tau} \cdot Q(\Bw, L(L'+\tau)).
$$

\end{lemma}

\begin{proof} Assume $\Bw=w_1\dots w_r$. We can write $v$ as $\frac{1}{\tau}
\sum_{i=1}^r a_i w_i$, where $|a_i| \le L$. An element $y$ in
$Q(\Bw, L) + Q(v, L')$ can be written as

$$ y= \sum_{i=1}^r x_i w_i + x v $$

\noindent where $|x_i| \le L$ and $|x| \le L'$. Substituting $v$,
we have

$$y= \sum_{i=1}^r x_i w_i + x \frac{1}{\tau} \sum_{i=1}^r a_i w_i
= \frac{1}{\tau} \sum_{i=1}^r w_i (\tau x_i + x a_i), $$

\noindent where $| \tau x_i + x a_i | \le \tau L + L'L$. This
concludes the proof. \end{proof}

\begin{proof}[of Theorem \ref{theo:inverse3}]
We begin by running the algorithm in the proof of Theorem
\ref{theo:inverse2} to locate a word $\Bw$ of length at most $d-1$
such that the set $\bigcup_{1 \leq \tau \leq k} \frac{1}{\tau}
\cdot Q(\Bw,k)$ covers all but at most $k^2$ elements of $\Bv$.
Set $\Bv^{[0]}$ be the word formed by removing the (at most $k^2$)
exceptional elements from $\Bv$ which do not lie in $\bigcup_{1
\leq \tau \leq k} \frac{1}{\tau} \cdot Q(\Bw,k)$.

By increasing the constant
 $k_0$ in the assumption of the theorem, we can assume, in all
 arguments below, that $k$ is sufficiently large, whenever needed.

 By \eqref{pmub}, \eqref{14}
\begin{equation}\label{pcd-recurse-0}
 \BBP_{\mu/4d} (\Bv^{[0]} \Bw^{k^2} ) \geq  \BBP_{\mu/4d} (\Bv \Bw^{k^2} )
\ge \BBP_{\mu/4d}  (\Bv ) \BBP_{\mu/4d} (\Bv \Bw^{k^2} ) \ge
k^{-d} \BBP_{\mu/4d} (\Bv \Bw^{k^2} ).
 \end{equation}

In the following, we assume that there is at least one non-zero
entry in $\Bw$, as otherwise the claim is trivial.

Now we perform an additional algorithm. Let $K = K(\mu,d, \ep) >
2$ be a large constant to be chosen later.

\begin{itemize}

\item Step 0.  Initialize $i = 0$ and  Set $Q_0:= Q(\Bw, k^2)$ and
$\Bv^{[0]}$ as above.

 \item Step
1. Count how many $v \in \Bv^{[i-1]}$ having torsion at least $K$
with respect to $2Q_{i-1}$. (We need to have the factor $2$ here
in order to apply Lemma \ref{tor}.) If this number is less than
$k^2$, halt the algorithm. Otherwise, move on to Step 2.

\item Step 2.   Locate a multiset $S$ of $k^2$ elements of
$\Bv^{[i-1]}$ with torsion at least $K$ with respect to
$2Q_{i-1}$. Applying \eqref{pmub-3a}, we can find an element $v
\in S$ such that
$$
\BBP_{\mu/4d}(\Bv^{[i-1]} \Bw^{k^2} W_1^{\tau_1^2} \ldots
W_{i-1}^{\tau_{i-1}^2}) \leq \BBP_{\mu/4d}(\Bv^{[i]} \Bw^{k^2}
W_1^{\tau_1^2} \ldots W_{i-1}^{\tau_{i-1}^2} v^{k^2} )$$ where
$\Bv^{[i]}$ is obtained from $\Bv^{[i-1]}$ by deleting $S$.

 Let
$\tau_i$ be the torsion of $v$ with respect to $2Q_{i-1}$. Since
every element of $\Bv^{[0]}$ has torsion at most $k$ with respect
to $Q_0$, $ K \le \tau_i \le k$. We then set $W_{i} := v$, $Q_i :=
Q_{i-1} + Q(W_i, \tau_i^2)$, increase $i$ to $i+1$ and return to
Step 1.

\end{itemize}

Consider a stage $i$ of the algorithm. From construction and
induction and \eqref{pcd-recurse-0}, we have a word $W_1 \ldots
W_i$ with
$$ \BBP_{\mu/4d}(\Bv^{[i]} \Bw^{k^2} W_1^{\tau_1 ^2} \ldots W_i^{\tau_i^2}) \geq \BBP( \Bv^{[0]} \Bw^{k^2} )
 \ge k^{-d} \BBP(\Bw^{k^2}).$$

On the other hand, by applying Lemma \ref{tor} iteratively, we
have

$$ \BBP_{\mu/4d}( \Bw^{k^2} W_1^{\tau_1^2} \ldots W_i^{\tau_i^2} ) \leq
\BBP_{\mu/4d}( \Bw^{k^2} ) \prod_{j=1}^i (C_{\mu} \tau_j^{-1}) .$$

\noindent It follows that $\prod_{j=1}^i (C_{\mu} \tau_j^{-1}) \ge
k^{-d}$, or equivalently $ \prod_{j=1}^i (C_{\mu}^{-1} \tau_j) \le
k^d$. Recall that $\tau_j \ge K$. Thus by setting $K$ sufficiently
large (compared to $C_{\mu}, d$ and $1/\ep$), we can guarantee
that

\begin{equation} \label{equT} \prod_{j=1}^i \tau_j \le k^{d+\ep/2d} \end{equation} where  $\ep$ is the
constant in the assumption of the theorem. It also follows that
the algorithm must terminate at some stage $D \le \log _K k^{d+
\ep/2d} \le (d+1) \log_K k $.

Let us take a look at the final set $Q_D$. Applying  Lemma
\ref{lemma:sumofQ} iteratively we have

$$Q_D \subset (\prod_{j=1}^D \frac{1}{\tau_j}) \cdot Q(\Bw, L_D) $$

\noindent where $L_0 := k^2$ and

\begin{equation} \label{equL} L_i := L_{i-1} (\tau_i + \tau_i^2)
\le (1+1/K) L_{i-1} \tau_i^2. \end{equation}

We now show that the GAP $Q:= \frac{1}{K!} \cdot (2K!)Q(\Bw,L_D)=
\frac{1}{K!} \cdot Q(\Bw, 2K!L_D)$
 satisfies the claims of the theorem.

\begin{itemize}

\item ({\it Rank}) We have $\rank (Q) = \rank (Q(\Bw, L_D)) =
\rank (Q_0)= r \le d-1$, as showed in the proof of the previous
theorem.

\item ({\it Volume}) We  have $\Vol (Q) = (2K!)^r \Vol (Q(\Bw,
L_D)) =O(\Vol (Q(\Bw, L_D)))$. On the other hand, by \eqref{equT}
and  \eqref{equL}

$$ \Vol
(Q(\Bw, L_D)) =  (2L_D +1 )^r \le (3L_D)^r = O( \Big(k^2
\prod_{j=1}^D (1+1/K) \tau_j^2 \Big)^r ) =O( \Big( k^{2+ 2(d+
\ep/2d) } (1+K)^D \Big)^r ). $$

By definition, $D \le \log_K k^{d+\ep/2d} < \log k$, given that
$K$ is sufficiently large compared to $d$. Thus  $(1+1/K)^D \le
\exp (D/K) \le k^{1/K } $ which implies that

$$ \Vol
(Q(\Bw, L_D)) =O( k^{r( 2 + 2 (d+ \ep/2d) + 1/K) }) =o(
k^{2(d^2-1) + \ep }) $$

\noindent provided that $r \le d-1$ and $K$ is sufficiently large
compared to $d$ and $1/\ep$. (The asymptotic notation here is used
under the assumption that $k \rightarrow \infty$.)

\item ({\it Number of exceptional elements}) At  each stage in the
second algorithm, we discard a set of $k^2$ elements, thus all but
$Dk^2 \le (d+1) k^2 \log_K k)$ elements of $\Bv^{[0]}$ have
torsion at most $K$ with respect to $ 2Q_D$. As $Q_D \subset
Q(\Bw, L_D)$ and $v \backslash v^{[0]} \le k^2$, it follows that
all but at most

$$(d+1) k^2 \log_K k + k^2 $$

\noindent elements of $\Bv$ have torsion at most $K$ with respect
to $2Q(\Bw, L_D) = Q(\Bw, 2L_D)$. By setting $K$ sufficiently
large compared to $d$ and $1/\ep$, we can guarantee that

$$(d+1) k^2 \log_K k + k^2  \le \ep k^2 \log k. $$

 To conclude, notice that  any element
with torsion at most $K$ with respect to $Q(\Bw, 2L_D)$ belongs to
$Q:=\frac{1}{K!} \cdot Q(\Bw, 2K!L_D)$. Thus, $Q$ contains all but
at most  $\ep k^2 \log k$ elements of $\Bv$.

\item ({\it Generators})  The generators of $\frac{1}{K!} \cdot
Q(\Bw, 2K!L_D)$ are $\frac{1}{K! \prod_{ j=1}^D \tau_j } w_i$, $1
\le i \le r $. Since $w_i \in \Bv$ and $\prod_{ j=1}^D \tau_j \le
k^{d+\ep/2d} =o(k^{d+\ep})$, the claim about generators follows.

\end{itemize} The proof is complete. \end{proof}

\section{The smallest singular value}\label{singsec}

In this section, we prove Theorem \ref{theo:singular1}, modulo two
key results, Theorem \ref{discrete-thm} and Corollary
\ref{cor:KKS2}), which will be proved in later sections.

Let $B
> 10$ be a large number (depending on $A$) to be chosen later.
Suppose that $\sigma_n(M^{\mu}_n) < n^{-B} $. This means that
there exists a unit vector $v$ such that
$$ \| M^{\mu}_n v \| < n^{-B}.$$

By rounding each coordinate $v$ to the nearest multiple of
$n^{-B-2}$, we can find a vector $\tilde v \in n^{-B-2} \cdot
\Z^n$ of magnitude $0.9 \le \|\tilde v \| \le 1.1$ such that
$$ \| M^{\mu} _n \tilde v \| \le 2n^{-B}.$$

Writing $w := n^{B+2} \tilde v$, we thus can find an integer
vector $w  \in \Z^n$ of magnitude $ 0.9 n^{B+2} \le\|w\| \le 1.1
n^{B+2} $ such that
$$ \| M^{\mu} _n w \| \le 2n^2.$$

Let $\Omega$ be the set of integer vectors $w  \in \Z^n$ of
magnitude $ 0.9 n^{B+2} \le\|w\| \le 1.1 n^{B+2} $. It suffices to
show the probability bound
$$ \P( \hbox{there is some } w \in \Omega \hbox{ such that } \| M^{\mu}_n w\| \le 2n^2 ) = O_{A,\mu}( n^{-A} ).$$
We now partition the elements $w = (w_1,\ldots,w_n)$ of $\Omega$ into three
sets:

\begin{itemize}

\item  We say  that $w$ is \emph{rich} if
$$ \BBP_{\mu}(w_1 \ldots w_n) \geq n^{-A-10} $$
and \emph{poor} otherwise. Let $\Omega_1$ be the set of poor
$w$'s.

\item A rich $w$ is \emph{singular} $w$ if fewer than $n^{0.2}$ of
its coordinates have absolute value $n^{B-10}$ or greater. Let
$\Omega_2$ be the set of rich and singular $w$'s.

\item A rich $w$ is \emph{non-singular} $w$, if at least $n^{0.2}$
of its coordinates have absolute value $n^{B-10}$ or greater. Let
$\Omega_3$ be the set of rich and non-singular $w$'s.

\end{itemize}

\noindent The desired estimate follows directly  from the
following lemmas and the union bound.

\begin{lemma}[Estimate for poor $w$]\label{lemma:Omega1}

$$\P( \hbox{there is some} \,\, w \in \Omega_1 \,\,\hbox{such
that} \,\,\|M^{\mu}_n w\| \le 2n^2 ) = o( n^{-A} ). $$

\end{lemma}

\begin{lemma}[Estimate for rich singular $w$]\label{lemma:Omega2}
$$\P( \hbox{there is some} \,\, w \in \Omega_2 \,\,\hbox{such
that} \,\,\|M^{\mu}_n w\| \le 2n^2 ) = o( n^{-A} ).
$$

\end{lemma}

\begin{lemma}[Estimate for rich non-singular $w$]\label{lemma:Omega3}
$$\P( \hbox{there is some} \,\, w \in \Omega_3 \,\,\hbox{such
that} \,\,\|M^{\mu}_n w\| \le 2n^2 ) = o( n^{-A} ).
$$

\end{lemma}

\begin{remark} Our arguments will show that the probabilities in
Lemmas \ref{lemma:Omega2} and \ref{lemma:Omega3} are exponentially
small. \end{remark}

The proofs of Lemmas \ref{lemma:Omega1} and \ref{lemma:Omega2} are
relatively simple and rely on well-known methods. We delay these
proofs to the end of this section and focus on the proof of Lemma
\ref{lemma:Omega3}, which is the heart of the matter, and which uses all the major tools discussed in previous
sections.

\begin{proof}[of Lemma \ref{lemma:Omega3}]  Informally, the strategy is to use the inverse Littlewood-Offord theorem
(Corollary \ref{cor:inverse3}) to place
the integers $w_1,\ldots,w_n$ in a progression, which we then discretize using Theorem \ref{discrete-thm}.  This allows us to
replace the event $\|M^{\mu}_n w\| \le 2n^2$ by the discretized event $M^{\mu,Y}_n = 0$ for a suitable $Y$, at which point
we apply Corollary \ref{cor:KKS2}.

We turn to the details.  Since $w$ is rich, we
see from Corollary \ref{cor:inverse3} that there exists a
symmetric GAP $Q$ of integers of rank at most $A'$ and volume at
most $n^{A'}$ which contains all but $\lfloor n^{0.1}\rfloor$ of
the integers $w_1,\ldots,w_n$, where $A'$ is a constant depending
on $\mu$ and $A$.  Also the generators of $Q$ are of the form $w_i/s$ for some
$1 \leq i \leq n$ and $1 \leq s \leq n^{A'}$.

Using the description of $Q$ and the fact that $w_1, \dots, w_n$
are polynomially bounded (in $n$), it is easy to derive that total
number of possible $Q$ is $n^{O_{A'}(1)}$. Next, by paying a factor of
$$ {n \choose \lfloor n^{0.1}\rfloor} \leq n^{\lfloor n^{0.1}\rfloor} = \exp(o(n))$$

\noindent we may assume that it is the last $\lfloor
n^{0.1}\rfloor$ integers $w_{m+1}, \ldots,w_n$ which possibly lie
outside $Q$, where we set $m := n-\lfloor n^{0.1} \rfloor$. As
each of the $w_i$ has absolute value at most $1.1 n^{B+2}$,  the
number of ways to fix these exceptional elements is at most $(2.2
n^{B+2} )^{n^{0.1}} = \exp(o(n))$. Overall, it costs a factor only
$\exp(o(n))$ to fix $Q$, the positions and values of the
exceptional elements of $w$.

Once we have fixed $w_{m+1}, \dots, w_n$,   we can then write
$$ M_n w = w_1 X^{\mu}_1 + \ldots + w_{m} X^{\mu}_{m} + Y,$$

\noindent where $Y$ is a random variable determined by $X^{\mu}_i$
and $w_i$, $m<i \le n$.  (In this proof we think of $X_i^{\mu}$ as
the column vectors of the matrix.) For any number $y$,  let $F_y$
be the event that there exists $w_1,\dots, w_m$ in $Q$, where at
least one of the $w_i$ has absolute value larger or equal
$n^{B-10}$, such that
$$ |w_1 X^{\mu}_1 + \ldots + w_{m} X^{\mu}_{m} + y| \le 2n^2 .$$

\noindent It suffices to prove that

$$\P( F_y) = o(n^{-A}) $$ for any $y$. Our argument will in fact
show that this probability is exponentially small.

We now apply Theorem \ref{discrete-thm} to the GAP $Q$ with $R_0
:= n^{B/2}$ and $S := n^{10}$ to find a scale $R = n^{B/2 +
O_A(1)}$ and symmetric GAPs $Q_{\operatorname{sparse}}$, $Q_{\operatorname{small}}$ of rank at
most $A'$ and volume at most $n^{A'}$ such that

\begin{itemize}

\item $ Q \subseteq Q_{\operatorname{sparse}} + Q_{\operatorname{small}}.$

\item $ Q_{\operatorname{small}} \subseteq [-n^{-10} R, n^{-10} R]$.

\item The elements of $n^{10} Q_{\operatorname{sparse}}$ are $n^{10}
R$-separated.

\end{itemize}

Since $Q$ (and hence $n^{10} Q$) contains $w_1,\ldots,w_m$, we can
therefore write
$$ w_j = w^{\operatorname{sparse}}_j + w^{\operatorname{small}}_j$$
for all $1 \leq j \leq m$, where $w^{\operatorname{sparse}}_j \in Q_{\operatorname{sparse}}$ and
$w^{\operatorname{small}}_j \in Q_{\operatorname{small}}$. In fact, this decomposition is
unique.

Suppose that the event $F_y$ holds.  Writing $X^{\mu}_i =
(\eta^{\mu}_{i,1}, \ldots, \eta^{\mu} _{i,n})$ (where
$\eta^{\mu}_{i,j} $ are, of course, i.i.d copies of $\eta^{\mu}$)
and $y = (y_1,\ldots,y_n)$, we have

$$ w_1 \eta^{\mu}_{i,1} + \ldots + w_m \eta^{\mu}_{i,m} = y_i + O(n^2).$$

\noindent for all $1 \leq i \leq n$.  Splitting the $w_j$ into
sparse and small components and estimating the small components
using the triangle inequality, we obtain

$$ w^{\operatorname{sparse}}_1 \eta^{\mu}_{i,1} + \ldots + w^{\operatorname{sparse}}_m \eta^{\mu}
_{i,m} = y_i + O(n^{-9} R)$$ for all $1 \leq i \leq n$.  Note that
the left-hand side lies in $mQ_{\operatorname{sparse}} \subset n^{10}
Q_{\operatorname{sparse}}$, which is known to be $n^{10} R$-separated.  Thus
there is a unique value for the right-hand side, call it $y'_i$,
which depends only on $y$ and $Q$
 such that
$$ w^{\operatorname{sparse}}_1 \eta_{i,1} + \ldots + w^{\operatorname{sparse}}_m \eta_{i,m} = y'_i.$$
The point is that we have now eliminated the $O()$ errors, and
have thus essentially converted the singular value problem to the zero
determinant problem.  Note also that since one of the
$w_1,\ldots,w_m$ is known to have magnitude at least $n^{B-10}$
(which will be much larger than $n^{10} R$ if $B$ is chosen large
depending on $A$), we see that at least  one of the
$w_1^{\operatorname{sparse}},\ldots,w_n^{\operatorname{sparse}}$ is non-zero.

Consider the random matrix $M'$ of order $m \times m+1$ whose
entries are i.i.d copies of $\eta^{\mu}$ and let $y' \in \R^{m+1}$
be the column vector $y' = (y'_1,\ldots,y'_{m+1})$. We conclude
that if the event $F_y$ holds, then there exists a non-zero vector
$w \in \R^m$ such that $M' w = y'$.  But from Corollary
\ref{cor:KKS2}, this holds with the desired probability
$$\exp(-\Omega(m+1)) =\exp(-\Omega(n)) =o(n^{-A})
$$ and we are done. \end{proof}

\begin{proof}[of Lemma \ref{lemma:Omega1}]
We use a conditioning argument, following \cite{Rud}. (An argument
of the same spirit was used by Koml\'os to prove the bound $O(n^{-1/2})$
for the singularity problem \cite{Bol}.)

Let $M$ be a matrix such that there is $w \in \Omega_1$ satisfying
$\|Mw\| \le 2n^2$. Since $M$ and its transpose have the same spectral norm, there is a vector $w'$ which has the
same norm as $w$ such that $\|w'M\| \le 2n^2$. Let $u=w'M$ and
$X_i$ be the row vectors of $M$. Then

$$u =\sum_{i=1}^n w_i' X_i $$

\noindent where $w_i'$ are the coordinates of $w'$.

Now we think of $M$ as a random matrix. By paying a factor of $n$,
we can assume that $w'_n$ has the largest absolute value among the
$w_i'$. We expose the first $n-1$ rows $X_1, \dots, X_{n-1}$ of
$M$. If there is  $w \in \Omega_1$ satisfying $\|Mw\| \le 2n^2$,
then there is a vector $y \in \Omega_1$, depending only on the
first $n-1$ rows such that

$$(\sum_{i=1}^{n-1} (X_i \cdot y)^2 )^{1/2} \le 2n^2. $$

\noindent Now consider the inner product $X_n\cdot y$. We can
write $X_n$ as

$$X_n =\frac{1}{w_n'} (u- \sum_{i=1}^{n-1} w_i' X_i ). $$

\noindent Thus,

$$|X_n \cdot y| = \frac{1}{\|w_n' \|} |u \cdot y  - \sum_{i=1}^{n-1} w_i' X_i \cdot y|.
$$

\noindent The right hand side, by the triangle inequality, is at
most

$$\frac{1}{\|w_n'\|} (\|u \| \|y \| + \|w'\| (\sum_{i=1}^{n-1} (X_i
\cdot y)^2 )^{1/2}). $$

\noindent By assumption $\|w_n' \| \ge n^{-1/2} \|w' \|$.
Furthermore, as $\|u \| \le 2n^2$, $\|u \| \|y \| \le 2n^2 \|y \|
\le 3n^2 \|w' \|$ as $\|w'\|=\|w \|$ and both $y$ and $w$ belong
to $\Omega_1$. (Any two vectors in $\Omega_1$ has roughly the same
length.) Finally $(\sum_{i=1}^{n-1} (X_i \cdot y)^2 )^{1/2} \le
2n^2$. Putting all these together, we have

$$|X_n \cdot y| \le 5n^{5/2}. $$

\noindent Recall that $y$ is fixed (after we expose the first
$n-1$ rows) and $X_n $ is a copy of $X^{\mu}$. The probability
that $|X^{\mu} \cdot y| \le 5n^{5/2}$ is at most $(10n^{5/2}+1)
\BBP_{\mu} (y)$. On the other hand, $y$ is poor, so $\BBP_{\mu}(y) \le
n^{-A-10}$. Thus, it follows that

$$\P( \hbox{there is some} \,\, w \in \Omega_1 \,\,\hbox{such
that} \,\,\|M^{\mu}_n w\| \le 2n^2 ) \le n^{-A-10} (10 n^{5/2}+1)
n = o( n^{-A} ),
$$

\noindent where the extra factor $n$ comes from the assumption
that $w_n'$ has the largest absolute value. This completes the
proof. \end{proof}

\begin{proof}[of Lemma \ref{lemma:Omega2}] We use an
argument from \cite{Lit}.  The key point will be that the set $\Omega_2$ of rich
non-singular vectors has sufficiently low entropy that one can proceed using the union bound.

A set $N$ of vectors on the
$n$-dimensional unit sphere $S_{n-1}$ is said to be an \emph{$\ep$-net} if for
any $x \in S_{n-1}$, there is $y\in N$ such that $\|x-y\| \le
\ep$.  A standard greedy argument shows

\begin{lemma} \label{lemma:net} For any $n$ and $\epsilon \le 1$, there exists an $\epsilon$-net
of cardinality at most $O(1/\eps)^n$.
\end{lemma}

Next, a simple concentration of measure argument shows

\begin{lemma} \label{lemma:smallnorm} For any fixed vector $y$ of magnitude between $0.9$ and $1.1$
$$\P (\|M^{\mu}_n y\| \le  n^{-2}) = \exp(-\Omega(n)). $$
\end{lemma}

It suffices to verify this statement for the case $|y|=1$. Notice
that

$$\|M^{\mu}_n y\|^2 = \sum_{i=1}^n (X_i \cdot y)^2 =\sum_{i=1}^n
Z_i  $$ where $Z_i = (X_i \cdot y)^2 $. The $Z_i$ are i.i.d random
variables with expectation $\mu$ and bounded variance. Thus
$\sum_{i=1}^n Z_i$ has mean $\Omega (n)$ and  the claimed bound
follows from Chernoff's large deviation inequality (see, e.g.,
\cite[Chapter 1]{TVbook}). (In fact, one can replace the $n^{-2}$
by $cn^{1/2}$ for some small constant $c$, but this refinement is
not necessary.)

For a vector $w \in \Omega_2$, let $w'$ be its normalization $w':=
w/\|w\|$. Thus, $w'$ is an unit vector with at most $n^{0.2}$
coordinates with absolute values larger or equal $n^{-10}$. Let
$\Omega_2'$ be the collection of those $w'$ with this property.

If $\| Mw\| \le 2n^2$ for some $w\in \Omega_2$, then $\|Mw'\| \le
3 n^{-B}$ , as $\|w\| \ge 0.9 n^{B+2}$. Thus, it suffices to give
an exponential bound on the event that there is $w'\in \Omega_2'$
such that $\| M^{\mu}_n w'\| \le 3n^{-B}$.

By paying a factor ${n \choose {n^{0.2}} }= \exp(o(n))$ in
probability, we can assume that the large coordinates (with
absolute value at least $n^{-10} $) are among the first $l :=
n^{0.2}$ coordinates. Consider an $n^{-3}$-net $N$ in $S_{l-1}$.
For each vector $y \in N$, let $y'$ be the $n$-dimensional vector
obtained from $y$ by letting the last  $n-l$ coordinates be zeros,
and let $N'$ be the set of all such vectors obtained. These
vectors have magnitude between $0.9$ and $1.1$, and from Lemma
\ref{lemma:net} we have $|N'| \leq O(n^3)^{l}$.

Now consider a rich singular vector $w' \in \Omega_2$ and let $w^{''}$ be the
$l$-dimensional vector formed by the first $l$ coordinates of this
vector. As the remaining coordinates are small $\|w^{''}\| =1
+O(n^{-9.5})$. There is a vector $y \in N$ such that

$$\|y - w^{''} \| \le n^{-3}+ O(n^{-9.5}) . $$

\noindent It follows that there is a vector $y' \in N'$ such that

$$\|y' - w'\| \le  n^{-3}+ O(n^{-9.5}) \le 2n^{-3}. $$

\noindent For any matrix $M$ of norm at most $n$

$$\|M w'\| \ge \|M y'\| - 2n^{-3} n = \|My'\| - 2n^{-2}. $$

\noindent It follows that if $\|M w'\| \le 3n^{-B}$ for some $B
\ge 2$, then $\| My'\| \le 5n^{-2}$. Now take $M=M^{\mu}_n$. For
each fixed $y'$, the probability that $\| My'\| \le 5n^{-2}$ is at
most $\exp(-\Omega(n))$, by Lemma \ref{lemma:smallnorm}. Furthermore,
the number of $y'$ is subexponential (at most $O(n^3)^{l} =
O(n)^{3 n^{.2}} = \exp(o(n))$). Thus the claim follows directly by
the union bound. \end{proof}

\section{Discretization of progressions}\label{discrete-sec}

The purpose of this section is to prove Theorem \ref{discrete-thm}.  The arguments here are elementary (based mostly on the pigeonhole principle and linear algebra, in particular Cramer's rule) and can be read independently of the rest of the paper.

We shall follow the informal strategy outlined in Section \ref{discrete-subsec}.
We begin with a preliminary observation, that basically asserts the intuitive fact that progressions do not contain large lacunary subsets.

\begin{lemma}\label{sgap} Let $P \subset \Z$ be a symmetric generalized arithmetic progression of rank $d$ and volume $V$, and let $x_1, \ldots, x_{d+1}$ be
non-zero elements of $P$.  Then there exist $1 \leq i < j \leq d+1$ such that
$$ C_d^{-1} V^{-1} |x_i| \leq |x_j| \leq C_d V |x_i|$$
for some constant $C_d > 0$ depending only on $d$.
\end{lemma}

\begin{proof} We may order $|x_{d+1}| \geq |x_d| \geq \ldots \geq |x_1|$.  If we write
$$ P = \{ m_1 v_1 + \ldots + m_d v_d: |m_i| \leq M_i \hbox{ for all } 1 \leq i \leq d \}$$
(so that $V = \Theta_d(M_1 \ldots M_d)$), then each of the $x_1,\ldots,x_{d+1}$ can be written as a linear combination of the $v_1,\ldots,v_d$.
Applying Cramer's rule, we conclude that there exists a non-trivial relation
$$ a_1 x_1 + \ldots + a_{d+1} x_{d+1} = 0$$
where $a_1,\ldots,a_{d+1} = O_d(V)$ are integers, not all zero.  If we let $j$ be the largest index such that $a_j$ is non-zero, then
$j > 1$ (since $x_1$ is non-zero) and we conclude in particular that
$$ |x_j| = O( |a_j x_j| ) = O_d( V |x_{j-1}| )$$
from which the claim follows.
\end{proof}

\begin{proof}[of Theorem \ref{discrete-thm}]
We can assume that $R_0$ is very large compared to $(SV)^{O_d(1)}$ since otherwise the claim is trivial (take $P_{\operatorname{sparse}} := P$
and $P_{\operatorname{small}} := \{0\}$).  We can also take $V \geq 2$.

Let $B = B_d$ be a large integer depending only on $d$ to be chosen later.
The first step is to subdivide the interval $[(SV)^{-B^{B+2}} R_0, (SV)^{B^{B+2}} R_0]$ into $\Theta(B)$ overlapping subintervals of the form $[(SV)^{-B^{B+1}} R, (SV)^{B^{B+1}} R]$, with every integer being contained in at most $O(1)$ of the subintervals.  From Lemma \ref{sgap} and the pigeonhole principle we see that at most $O_d(1)$ of the intervals can contain an element of $(SV)^{B^B} P$ (which has volume $O((SV)^{O_d(B^B)})$.  If we let $B$ be sufficiently large, we can thus find an interval
$[(SV)^{-B^{B+1}} R, (SV)^{B^{B+1}} R]$ which is disjoint from $(SV)^{B^B} P$.  Since $P$ is symmetric, this means that every $x \in (SV)^{B^B} P$ is either larger than $(SV)^{B^{B+1}} R$ in magnitude, or smaller than $(SV)^{-B^{B+1}} R$ in magnitude.

Having located a good scale $R$ to discretize, we now split $P$ into small ($\ll R$) and sparse ($\gg R$-separated) components.
We write $P$ explicitly as
$$ P = \{ m_1 v_1 + \ldots + m_d v_d: |m_i| \leq M_i \hbox{ for all } 1 \leq i \leq d \}$$
so that $V = \Theta_d(M_1 \ldots M_d)$ and more generally
$$ kP = \{ m_1 v_1 + \ldots + m_d v_d: |m_i| \leq kM_i \hbox{ for all } 1 \leq i \leq d \}$$
for any $k \geq 1$.  For any $1 \leq s \leq B$, let $A_s \subset \Z^d$ denote the set
$$ A_s := \{ (m_1,\ldots,m_d): |m_i| \leq V^{B^s} M_i \hbox{ for all } 1 \leq i \leq d; |m_1 v_1 + \ldots + m_d v_d| \leq (SV)^{-B^{B+1}} R \}.$$
Roughly speaking, this space corresponds to the kernel of $\Phi$
as discussed in Section \ref{discrete-subsec}; the additional parameter $s$ is a technicality needed to compensate for the fact that boxes, unlike vector spaces, are not quite closed under dilations.
We now view $A_s$ as a subset of the Euclidean space $\R^d$.  As such it spans a vector space $X_s \subset \R^d$.  Clearly
$$ X_1 \subseteq X_2 \subseteq \ldots \subseteq X_B$$
so if $B$ is large enough, then by the pigeonhole principle (applied to the dimensions of these vector spaces) we can find
$1 \leq s < B$ such that we have the stabilization property $X_s = X_{s+1}$.  Let the dimension of this space be $r$, thus $0 \leq r \leq d$.

There are two cases, depending on whether $r=d$ or $r<d$.  Suppose first that $r=d$ (so the kernel has maximal dimension).
Then by definition of $A_s$ we
have $d$ ``equations'' in $d$ unknowns,
$$ m^{(j)}_1 v_1 + \ldots + m^{(j)}_d v_d = O( (SV)^{-B^{B+1}} R ) \hbox{ for all } 1 \leq j \leq d,$$
where $m^{(j)}_i = O( M_i V^{B^s} )$ and the vectors $(m^{(j)}_1,\ldots,m^{(j)}_d) \in A_s$ are linearly independent as $j$ varies.
Using Cramer's rule we conclude that
$$ v_i = O_d( (SV)^{O_d(B^s)} (SV)^{-B^{B+1}} R ) \hbox{ for all } 1 \leq j \leq d$$
since all the determinants and minors which arise from Cramer's rule are integers that vary from $1$ to $O_d(V^{O_d(B)})$ in
magnitude.  Since $M_i = O(V)$ for all $i$, we conclude that $x = O_d( V^{O_d(B^s)} (SV)^{-B^{B+1}} R )$ for all $x \in P$, which by construction of $R$ (and the fact that $s<B$) shows that
$P \subset [-(SV)^{-B^{B+1}} R, (SV)^{-B^{B+1}} R]$ (if $B$ is sufficiently large).  Thus in this case we can take $P_{\operatorname{small}} = P$ and $P_{\operatorname{sparse}} = \{0\}$.

Now we consider the case when $r < d$ (so the kernel is proper).  In this case we can write $X_s$ as a graph of some linear transformation $T: \R^r \rightarrow  \R^{d-r}$: after
permutation of the coordinates, we have
$$ X_s = \{ (x, Tx) \in \R^r \times \R^{d-r}: x \in \R^r \}.$$
The coefficients of $T$ form an $r \times d-r$ matrix, which can be computed by Cramer's rule to be rational numbers with numerator and denominator
$O_d((SV)^{O_d(B^s)})$; this follows from $X_s$ being spanned by $A_s$, and on the integrality and size bounds on the coefficients of elements of $A_s$.

Let $m \in A_s$ be arbitrary.  Since $A_s$ is also contained in $X_s$, we can
write $m = (m_{[1,r]}, T m_{[1,r]})$ for some $m_{[1,r]} \in \Z^r$ with magnitude
$O_d((SV)^{O_d(B^s)})$.  By definition of $A_s$,
we conclude that
$$ \langle m_r, v_{[1,r]} \rangle_{\R^r} + \langle T m_r, v_{[r+1,d]} \rangle_{\R^{d-r}} = O( (SV)^{-B^{B+1}} R )$$
where $v_{[1,r]} := (v_1,\ldots,v_r)$, $v_{[r+1,d]} := (v_{r+1},\ldots,v_d)$, and the inner products on $\R^r$ and $\R^{d-r}$ are the
standard ones.  Thus
$$ \langle m_r, v_{[1,r]} + T^* v_{[r+1,d]} \rangle_{\R^{r}} = O( (SV)^{-B^{B+1}} R )$$
where $T^*: \R^{d-r} \rightarrow \R^r$ be the adjoint linear
transformation to $T$.  Now since $A$ spans $X$, we see that the
$m_{[1,r]}$ will linearly span $\R^r$ as we vary over all elements
$m$ of $A$. Thus by Cramer's rule we conclude that
\begin{equation}\label{vtw}
 v_{[1,r]} + T^* v_{[r+1,d]} = O_d( V^{O_d(B^s)} (SV)^{-B^{B+1}} R ).
\end{equation}

Write $(w_1,\ldots,w_r) := T^* v_{[r+1,d]}$, thus $w_1,\ldots,w_r$ are rational numbers.
We then construct the symmetric generalized arithmetic progressions $P_{\operatorname{small}}$ and $P_{\operatorname{sparse}}$ explicitly as
$$ P_{\operatorname{sparse}} := \{ m_1 w_1 + \ldots + m_r w_r + m_{r+1} v_{r+1} + \ldots + m_d v_d: |m_i| \leq M_i \hbox{ for all } 1 \leq i \leq d \}$$
and
$$ P_{\operatorname{small}} := \{ m_1 (v_1+w_1) + \ldots + m_r (v_r+w_r): |m_i| \leq M_i \hbox{ for all } 1 \leq i \leq d \}.$$
It is clear from construction that $P \subseteq P_{\operatorname{sparse}} + P_{\operatorname{small}}$, and that $P_{\operatorname{sparse}}$ and $P_{\operatorname{small}}$ have
rank at most $d$ and volume at most $V$.  Now from \eqref{vtw} we have
$$ v_i + w_i = O_d( (SV)^{O_d(B^s)} (SV)^{-B^{B+1}} R )$$
and hence for any $x \in P_{\operatorname{small}}$ we have
$$ x = O_d( (SV)^{O_d(B^s)} (SV)^{-B^{B+1}} R ).$$
By choosing $B$ large enough we conclude
$$ |x| \leq R/S$$
which gives the desired smallness bound on $P_{\operatorname{small}}$.

The only remaining task is to show $SP_{\operatorname{sparse}}$ is sparse.  It suffices to show that $SP_{\operatorname{sparse}}-SP_{\operatorname{sparse}}$ has no non-zero intersection with
$[-RS, RS]$.  Suppose for contradiction that this failed. Then we can find $m_1,\ldots,m_d$ with $|m_i| \leq 2SM_i$ for all $i$ and
$$ 0 < m_1 w_1 + \ldots + m_r w_r + m_{r+1} v_{r+1} + \ldots + m_d v_d < RS.$$
Let $Q$ be the least common denominator of all the coefficients of $T^*$, then $Q = O_d( (SV)^{O_d(B^s)} )$.  Multiplying the above equation by $Q$,
we obtain
$$ 0 < m_1 Q w_1 + \ldots + m_r Q w_r + m_{r+1} Q v_{r+1} + \ldots + m_d Q v_d < O( RS V^{O_d(B^s)} ) < (SV)^{B^{B+1}} R.$$
Since $(w_1,\ldots,w_r) = T^* v_{[r+1,r+d]}$, the expression between the inequality signs is an integer linear combination of
$v_{r+1},\ldots,v_d$, with all coefficients of size $O_d( (S V)^{O_d(B^s)} )$, say
$$ m_1 Q w_1 + \ldots + m_r Q w_r + m_{r+1} Q v_{r+1} + \ldots + m_d Q v_d = a_{r+1} v_{r+1} + \ldots + a_d v_d.$$
In particular we see that this expression lies in $(SV)^{B^B} P$ (again
taking $B$ to be sufficiently large).  Thus by construction of $R$, we can improve the upper bound of $(SV)^{B^{B+1}}R$ to $(SV)^{-B^{B+1}} R$:
\begin{equation}\label{arv}
0 < a_{r+1} v_{r+1} + \ldots + a_d v_d  < (SV)^{-B^{B+1}} R.
\end{equation}
Taking $B$ to be large, this implies that $(0,\ldots,0,a_{r+1},\ldots,a_d)$ lies in $X_{s+1}$, which equals $X_s$.  But $X_s$ was a graph
from $\R^r$ to $\R^{d-r}$, and thus $a_{r+1} = \ldots = a_d = 0$, which contradicts \eqref{arv}.  This establishes the sparseness.
\end{proof}

\section{Proof of Theorem \ref{theo:KKS2}}\label{kks2-sec}

Let $Y=\{y_1, \dots, y_l\}$ be a set of $l$ independent vectors in
$\R^n$. Let us recall that  $M_n^{\mu, Y}$ denote the random
matrix with row vectors $X_1^{\mu}, \dots, X^{\mu}_{n-l}, y_1,
\dots, y_l$, where $X_i^{\mu}$ are i.i.d copies of $X^{\mu}
=(\eta^{\mu}_1 \dots, \eta^{\mu}_n)$.

Define $\delta(\mu) := \max\{1-\mu, \mu/2\}$. It is easy to show
that for any subspace $V$ of dimension $d$
\begin{equation}\label{equ:deltamu} \P(X^{\mu} \in V) \le
\delta(\mu)^{d-n}. \end{equation}

In the following, we are going to use $N$ to denote the quantity
$(1/\delta(\mu))^n$. As $0 < \mu \le 1$, $\delta(\mu) >0$ and thus
$N$ is exponentially large in $n$.  Thus it will suffice to show that
$$\P(M_n^{\mu,Y} \,\,\hbox{singular}\,\,) \le N^{-\eps+ o(1)}$$
for some $\eps=\eps(\mu,l) > 0$, where the $o(1)$ term is allowed to depend on $\mu$, $l$, and $\eps$.
  We may assume that $n$ is large depending on $\mu$ and $l$
since the claim is trivial otherwise.

Notice that if $M_n^{\mu,Y}$ is singular, then the row vectors
span a proper subspace $V$.  To prove the theorem, it suffices to
show that for any sufficiently small positive constant $\eps$

$$ \sum_{V, V \hbox{proper subspace}}   \P(
X_1^{\mu},\ldots,X_{n-l}^{\mu}, y_1, \dots, y_l \,\, \hbox{span}
\,\, V) \le N^{-\eps+o(1)}.
$$

Arguing as in \cite[Lemma 5.1]{TV1}, we can restrict ourselves
to hyperplanes. Thus, it is enough to prove
$$ \sum_{V, V \hbox{hyperlane}}   \P(
X_1^{\mu},\ldots,X_{n-l}^{\mu}, y_1, \dots, y_l \,\, \hbox{span}
\,\, V) \le N^{-\eps + o(1)}.
$$

 Clearly, we may restrict our attention to those
hyperplanes $V$ which are spanned by their intersection with
$\{-1,0,1\}^n$, together with $y_1,\ldots,y_l$. Let us
call such hyperplanes \emph{non-trivial}.
Furthermore, we call a hyperplane $H$ \emph{degenerate} if there
is a vector $v$ orthogonal to $H$ and at most $\log\log n$
coordinates of $v$ are non-zero.  Following
\cite[Lemma 5.3]{TV1}, it is easy to see that the number of degenerate
non-trivial hyperplanes is at most $N^{o(1)}$. Thus, their
contribution in the sum is at most
$$N^{o(1)} \delta(\mu)^{n-l} = N^{-1+o(1)}$$
\noindent which is acceptable. Therefore, from now on we
can assume that $V$ is non-degenerate.

 For each non-trivial hyperplane $V$, define the \emph{discrete
codimension} $d(V)$ of $V$ to be the unique integer multiple of
$1/n$ such that
\begin{equation}\label{codimension-def}
 N^{-\frac{d(V)}{n}-\frac{1}{n^2}} < \P( X^{\mu} \in V) \leq
N^{-\frac{d(V)}{n}}.
\end{equation}
Thus $d(V)$ is large when $V$ contains few elements from
$\{-1,0,1\}^n$, and conversely.

Let $B_V$ denote the event that $X_1^{\mu}, \dots, X_{n-l}^{\mu},
y_1, \dots, y_l$ span $V$.  We denote by $\Omega_d$ the set of all
non-degenerate, non-trivial hyperplanes with discrete codimension
$d$. It is simple to see that $ 1 \leq d(V) \leq n^2$ for all
non-trivial $V$. In particular, there are $n^2 = N^{o(1)}$
possible values of $d$, so to prove our theorem it suffices to
show that
\begin{equation}\label{d-dim}
 \sum_{V \in \Omega_d} \P( B_V)  \leq N^{-\eps +o(1)}
\end{equation}
for all $1 \leq d \leq n^2$.

We first handle the (simpler) case when $d$ is large. Note that if
$X_1^{\mu},\ldots,X_{n-l}^{\mu}, y_1, \dots, y_l$ span $V$, then
some subset of $n-l-1$ vectors $X_i$ together with the $y_j$'s
already span $V$ (since the $y_j$'s are independent). By symmetry,
we have

\begin{align*}
\sum_ {V \in \Omega_d} \P( B_V) &\leq (n-l) \sum_{V \in \Omega_d}
\P( X^{\mu}_1,\ldots,X^{\mu}_{n-l-1}, y_1, \dots, y_l  \hbox{ span
} V) \P(X_{n-l}^{\mu} \in V)
\\ &\le nN^{-\frac{d}{n}}\sum_{V \in \Omega_d}  \P(
X_1^{\mu},\ldots,X_{n-l-1}^{\mu}, y_1, \dots, y_l  \hbox{ span } V) \\
&\le nN^{-\frac{d}{n}} = N^{-\frac{d}{n} + o(1)}
\end{align*}

This disposes of the case when $d \geq \eps n$.  It remains to
verify the following lemma.

\begin{lemma} \label{mainbound} For all sufficiently small positive constant $\eps$, the following holds.
 If $d$ is any integer multiple of $1/n$ such that
\begin{equation}\label{d-bound}
1 \leq d \leq (\eps - o(1))n
\end{equation}
then we have
$$\sum_{V \in \Omega_d}
\P( B_V) \le  N^{-\eps+o(1)}.$$
\end{lemma}

\begin{proof}
For $0 < \mu \le 1$ we define the quantity $0 < \mu^{\ast} \leq 1/8$ as follows. If $\mu=1$
then $\mu^{\ast}:=1/16$. If $1/2 \le \mu <1$, then $\mu^{\ast}:=
(1-\mu)/4$. If $0 < \mu < 1/2$, then $\mu^{\ast}:= \mu/4$. We will need
the following inequality, which is a generalization of \cite[Lemma 6.2]{TV1}.

\begin{lemma}\label{hl}  Let $V$ be a non-degenerate
non-trivial hyperplane.  Then we have
$$ \P( X^{\mu} \in V ) \leq (\frac{1}{2} + o(1)) \P( X^{\mu^{\ast}} \in V ).$$
\end{lemma}

The proof of Lemma \ref{hl} relies on some Fourier-analytic ideas of Hal\'asz \cite{Hal} (see also \cite{KKS}, \cite{TV1}, \cite{TV2}) and is deferred till the end of the section.
Assuming it for now, we continue the proof of Lemma \ref{mainbound}.

Let us set $\gamma := \frac{1}{2}$; this is not the optimal value of this parameter, but will
suffice for this argument.

Let  $A_V$ be the event  that $
X^{\mu^{\ast}}_1,\ldots,X^{\,\mu^{\ast}}_{(1-\gamma) n}, \overline
X^{\mu}_1, \ldots, \overline X^{\mu}_{(\gamma-\eps) n}$ are
linearly independent in $V$, where $X^{\mu^{\ast}}_i$'s are i.i.d
copies of $X^{\mu^{\ast}}$ and $\overline X^{\mu}_j$'s are i.i.d
copies of $X^{\mu}$.

\begin{lemma}\label{lemma:AV} $$  \P(A_V) \ge
N^{(1-\gamma) -(1-\eps)d +o(1)}. $$ \end{lemma}

\begin{proof}
Notice that the right hand side on the bound in Lemma
\ref{lemma:AV} is the probability of the event $A'_V$ that $
X^{\mu^{\ast}}_1, \ldots, X^{\mu^{\ast}}_{(1-\gamma)n}, \overline
X^{\mu}_1, \ldots, \overline X^{\mu}_{(\gamma - \eps)n}$ belong to
$ V$. Thus, by Bayes' identity it is sufficient to show that
$$\P(A_V|A'_V) = N^{o(1)}. $$

 From \eqref{codimension-def} we
have
\begin{equation}\label{V-1}
 \P( X^{\mu} \in V ) = (1 + O(1/n)) \delta(\mu)^d
 \end{equation}
and hence by Lemma \ref{hl}
\begin{equation}\label{V-16}
 \P( X^{\mu^{\ast}} \in V ) \geq (2 + O(1/n)) \delta(\mu)^d.
 \end{equation}
On the other hand, by \eqref{equ:deltamu}

$$ \P( X^{\mu^{\ast}} \in W ) \leq (1-\mu^{\ast})^{n - \dim(W)}$$
for any subspace $W$.  By Bayes' identity we thus have the
conditional probability bound
$$ \P( X^{\mu^{\ast}} \in W | X^{(\mu^{\ast})} \in V) \leq
(2 +O(1/n))^{-1} \delta(\mu)^{-d}  (1-\mu^{\ast})^{n - \dim(W)}
\le \delta(\mu)^{-d}  (1-\mu^{\ast})^{n - \dim(W)}.$$

\noindent  When $\dim(W) \leq (1-\gamma)n$ the bound is less than
one when $\eps$ is sufficiently small, thanks to the bound on $d$ and the choice
$\gamma = \frac{1}{2}$.

Let $E_k$ be the event that $X^{\mu^{\ast}}_1, \ldots,
X^{\mu^{\ast}}_k$ are linearly independent. The above estimates
imply that

$$\P(E_{k+1} | E_k \wedge A'_V)
 \geq
1 -  \delta(\mu)^{-d} (1-\mu^{\ast})^{n - k}. $$

\noindent for all $0 \leq k \leq (1-\gamma)n$. Applying Bayes'
identity repeatedly  we thus obtain
$$  \P( E_{(1-\gamma)n} | A'_V ) \geq
N^{-o(1)}.$$

To complete the proof, observe that since
$$ \P( X^{\mu} \in W ) \leq \delta(\mu)^{n - \dim(W)}$$
for any subspace $W$, and hence by \eqref{V-1}
$$ \P( X ^{\mu} \in W | X^{\mu} \in V) \leq (1 + O(1/n)) \delta(\mu)^{-d}  \delta(\mu)^{n-\dim(W)}.$$

Let us assume $E_{(1-\gamma) n}$ and denote by $W$ the $(1-\gamma)
n$-dimensional subspace spanned by $X^{\mu^{\ast}}_1, \ldots,
X^{\mu^{\ast}}_{(1-\gamma)n}$. Let $U_k$ denote the event that $
\overline X^{\mu}_1, \ldots, \overline X^{\mu}_{k}, W \hbox{ are
liearly independent}$. We have

$$p_k= \P(U_{k+1}|U_k \wedge A'_V) \ge  1 - (1 + O(1/n)) \delta(\mu)^{-d} \delta(\mu)^{n-k-(1-\gamma)n} \ge
 1 - \frac{1}{100} \delta(\mu)^{-(\gamma-\eps)n+k} $$

\noindent  for all $0 \leq k < (\gamma - \eps)n$, thanks to
\eqref{d-bound}. Thus by Bayes' identity we obtain

$$\P(A_V|A'_V) \ge N^{o(1)}  \prod_{0 \leq k < (\gamma - \eps)n} p_k =
N^{o(1)} $$  as desired. \end{proof}

Now we continue the proof of the theorem.  Fix $V \in \Omega_d$.
Since $A_V$ and $B_V$ are independent, we have, by Lemma
\ref{lemma:AV}  that

$$\P(B_V) =\frac{\P(A_V \wedge B_V)}{\P(A_V)} \le N^{- (1-\gamma) +
(1-\eps)d + o(1)} \P(A_V \wedge B_V). $$

Consider a set
$$X^{\mu^{\ast}}_1,\ldots,X^{\mu^{\ast}}_{(1-\gamma) n}, \overline X^{\mu}_1, \ldots,
\overline X^{\mu}_{(\gamma-\eps) n},
X^{\mu}_1,\ldots,X^{\mu}_{n-l}$$ of vectors satisfying $A_V \wedge
B_V$. Then there exists $\eps n -l- 1$ vectors $X^{\mu}_{j_1},
\ldots, X^{\mu}_{j_{\eps n-l-1}}$ inside
$X^{\mu}_1,\ldots,X^{\mu}_{n-l}$ which, together with
$$X^{\mu^{\ast}}_1,\ldots,X^{\mu^{\ast}}_{(1-\gamma) n},
\overline X^{\mu}_1, \ldots, \overline X^{\mu}_{(\gamma-\eps) n},
y_1, \dots, y_l $$ span $V$. Since the number of possible indices
$j_1,\ldots, j_{\eps n-l-1}$ is ${n-l \choose \eps n - l-1} =
2^{(h(\eps) + o(1))n}$ (with $h$ being the entropy function), by
conceding a factor of

 $$2^{(h(\eps ) +o(1))n} = N^{ah(\eps ) +o(1)}, $$

 \noindent where $a=\log_{1/\delta(\mu)} 2$,
we can assume that $j_i =i$ for all relevant $i$. Let $C_V$ be the
event that $$X^{\mu^{\ast}}_1,\ldots,X^{\mu^{\ast}}_{(1-\gamma)
n}, \overline X^{\mu}_1, \ldots, \overline X^{\mu}_{(\gamma-\eps)
n}, X^{\mu}_1, \ldots, X^{\mu}_{\eps n-l-1}, y_1, \dots, y_l
\hbox{ span } V. $$ Then we have

$$\P(B_V) \le  N^{-(1-\gamma) +
(1-\eps)d + ah(\eps) + o(1)} \P\Big(C_V \wedge (X^{\mu}_{\eps
n},\ldots,X^{\mu}_{n -l} \hbox{ in } V ) \Big). $$

\noindent On the other hand,  $C_V$ and the event $(X_{\eps
n},\ldots,X_n \hbox{ in } V )$ are independent, so

$$\P\Big(C_V \wedge (X^{\mu} _{\eps n},\ldots,X^{\mu} _{n-l} \hbox{ in } V ) \Big)
=\P(C_V) \P( X^{\mu} \in V )^{(1-\eps)n+1-l}. $$

\noindent Putting the last two estimates together we obtain

\begin{align*}  \P(B_V) &\leq N^{-(1-\gamma) + (1-\eps)d + ah(\eps) + o(1)}
N^{-((1-\eps) n + 1-l)d/n} \P(C_V) \\ & =N^{-(1-\gamma) +ah(\eps)
+(l-1)\eps+o(1)} \P(C_V). \end{align*}

\noindent Since any set of vectors can only span a single space
$V$, we have $\sum_{V \in \Omega_d} \P(C_V) \le 1$. Thus, by
summing over $\Omega_d$, we have

$$\sum_{V \in \Omega_d} \P(B_V) \le N^{-(1-\gamma) +ah(\eps) +(l-1) \eps+o(1)}. $$

\noindent With the choice $\gamma = \frac{1}{2}$, we obtain a bound of
$N^{-\eps + o(1)}$ as desired, by choosing $\eps$ sufficiently small.
This provides the desired bound in Lemma
\ref{mainbound}.
\end{proof}

\subsection{Proof of Lemma \ref{hl}}

To conclude, we prove Lemma \ref{hl}. Let $v=(a_1, \dots, a_n)$ be
the normal vector of $V$ and define $$F_{\mu} (\xi) :=
\prod_{i=1}^n ((1-\mu)+ \mu \cos 2\pi a_i \xi). $$
{}From Fourier analysis we have (cf. \cite{TV1})
$$\P(X^{\mu} \in V) = \P(X^{\mu} \cdot v =0) =\int_{0}^1 F_{\mu}(\xi) d \xi. $$

The proof of Lemma \ref{hl}  is based on the following technical
lemma.

\begin{lemma} \label{3estimates}Let $\mu_1$ and $\mu_2$ be a positive numbers at most $1/2$
such that
the following two properties hold for  for  any $\xi, \xi' \in
[0,1]$.

\begin{equation}\label{fag}
 F_{\mu_1}(\xi) \leq F_{\mu_2}(\xi)^4
 \end{equation}
and
\begin{equation}\label{ffg}
F_{\mu_1}(\xi) F_{\mu_1}(\xi') \leq F_{\mu_2}(\xi + \xi')^2.
\end{equation}

\noindent Furthermore,
\begin{equation}\label{G-int}
 \int_0^1 F_{\mu_1}(\xi)\ d\xi = o(1). \end{equation} Then

 \begin{equation} \label{G-market}  \int_0^1 F_{\mu_1}(\xi)\ d\xi  \le (1/2+o(1)) \int_0^1 F_{\mu_2}(\xi)\
 d\xi. \end{equation}
\end{lemma}

\begin{proof} Notice that since $\mu_1, \mu_2 \le 1/2$, $F_{\mu_1}(\xi)$ and $F_{\mu_2}(\xi)$
are positive for any $\xi$.  From \eqref{ffg} we have the sumset
inclusion
$$
\{ \xi \in [0,1]: F_{\mu_1}(\xi) > \alpha \} + \{ \xi \in [0,1]:
F_{\mu_1}(\xi)
> \alpha \} \subseteq \{ \xi \in [0,1]: F_{\mu_2}(\xi) > \alpha \}$$ for
any $\alpha > 0$.  Taking measures of both sides and applying the
Mann-Kneser-Macbeath ``$\alpha+\beta$ inequality'' $|A+B| \geq
\min(|A|+|B|, 1)$ (see \cite{Mac}), we obtain
$$ \min( 2|\{ \xi \in [0,1]: F_{\mu_1}(\xi) > \alpha \}|, 1)
\leq  |\{ \xi \in [0,1]: F_{\mu_2}(\xi) > \alpha \}|.$$ But from
\eqref{G-int} we see that $|\{ \xi \in [0,1]: F_{\mu_2}(\xi) >
\alpha \}|$ is strictly less than 1 if $\alpha > o(1)$.  Thus we
conclude that
$$ |\{ \xi \in [0,1]: F_{\mu_1}(\xi) > \alpha \}| \leq \frac{1}{2} |\{ \xi \in
[0,1]: F_{\mu_2}(\xi) > \alpha \}|$$ when $\alpha > o(1)$.
Integrating this in $\alpha$, we obtain
$$
\int_{[0,1]: F_{\mu_1}(\xi) > o(1)} F_{\mu_1} (\xi)\ d\xi \leq
\frac{1}{2} \int_0^1 F_{\mu_2}(\xi)\ d\xi.$$ On the other hand,
from \eqref{fag} we see that when $F_{\mu_1}(\xi) \leq o(1)$, then
$F_{\mu_1}(\xi) = o( F_{\mu_1}(\xi)^{1/4}) =o(F_{\mu_2}(\xi))$,
and thus
$$
\int_{[0,1]: F_{\mu_1}(\xi) \leq o(1)} F_{\mu_1}(\xi)\ d\xi \leq
o(1) \int_0^1 F_{\mu_2}\ d\xi.$$ Adding these two inequalities we
obtain \eqref{G-market} as desired. \end{proof}

By Lemma \ref{lemma:properties}

$$\P(X^{\mu} \cdot v=0) \le \BBP_{\mu}(\Bv) \le \BBP_{\mu/4} (\Bv) =
\int_{0}^1 F_{\mu/4} (\xi) d \xi. $$

 It suffices to show that the conditions of Lemma \ref{3estimates} hold with
 $\mu_1= \mu/4$ and  $\mu_2= \mu^{\ast} =\mu/16$. The last estimate $ \int_0^1
F_{\mu_1}(\xi)\ d\xi \leq o(1)$ is a simple corollary of the fact
that at least $\log \log n$ among the $a_i$ are non-zero (instead
of $\log \log n$, one can use any function tending to infinity
with $n$), so we only need to verify the other two. Inequality
\eqref{fag} follows from the fact that $\mu_2=\mu_1/4$ and the
proof of the fourth property of Lemma \ref{lemma:properties}.

To verify \eqref{ffg}, we suffices to show that for any $\mu' \le
1/2$ and  any $\theta, \theta'$

$$((1-\mu')+ \mu' \cos \theta)((1-\mu')+ \mu' \cos \theta') \le
((1-\mu'/4)+ \frac{\mu'}{4} \cos (\theta + \theta')^2. $$

\noindent The left hand side is bounded from above by $((1-\mu')+
\mu' \cos \frac{\theta+\theta'}{2})^2 $, due to convexity. Thus, it
remains to show that

$$(1-\mu')+ \mu' \cos \frac{\theta+\theta'}{2} \le (1-\frac{\mu'}{4})+
\frac{\mu'}{4} \cos (\theta + \theta') $$

\noindent since both expressions are positive for $\mu' < 1/2$. By
defining $x:= \cos \frac{\theta+ \theta'}{2}$, the last inequality becomes

$$(1-\mu') + \mu' x \le (1-\frac{\mu'}{4}) +\frac{\mu'}{4} (2x^2-1) $$
which trivially holds. This completes the proof of Lemma \ref{hl}.

\end{document}